%% file: manuscript.tex
\documentclass[11pt,onecolumn,journal]{IEEEtran} 
\usepackage[cmex10]{amsmath}
\usepackage{amssymb,amsthm}			
\usepackage{bm,xspace,fancyhdr}					
\usepackage{graphicx}		
\usepackage[tight,footnotesize]{subfigure}
\usepackage{epstopdf}                
\usepackage{cite,url}                
\usepackage{pdfsync}				
\usepackage{setspace}
\graphicspath{{./figs/}}

\input{def.tex}

\begin{document}



\title{Analog Sparse Approximation with Applications to Compressed Sensing}%
\author{Adam~S.~Charles,~\IEEEmembership{Student Member,~IEEE},
        Pierre~Garrigues,
	and~Christopher~J.~Rozell$^{*}$,~\IEEEmembership{Member,~IEEE}
\thanks{Manuscript received November 4, 2011.  This work was supported in part by NSF grant CCF-0905346.}%
\thanks{* Corresponding author. ASC and CJR are with the School of Electrical and Computer Engineering, Georgia Institute of Technology, Atlanta, GA, 30332-0250 USA (e-mail: \{acharles6,crozell\}@gatech.edu).  PG is with IQ Engines, Berkeley, CA, 94704 (e-mail: pierre@iqengines.com). Preliminary versions of portions of this work were presented in~\cite{ROZ:2010b}.  The authors are grateful to Bruno Olshausen, Justin Romberg, Paul Hasler and Sam Shapero for valuable discussions related to this work, and to Marijn Brummer, Emory University and Children's Healthcare of Atlanta for providing MRI data. }
}

\maketitle

\begin{abstract}
\input{abstract}

\end{abstract}


\input{intro}
\input{background}
\input{BPDNperformance}

\input{other_costs}

\input{conclusions}

\bibliographystyle{IEEEtran}           

\bibliography{IEEEabrv,Adam_Bibliography2,RWL1,BPDNBIB} 

\input{apx_logbarrier}

\end{document}

%% file: def.tex
\newcommand{\bmath}[1]	{\ensuremath{\bm{#1}}\xspace}

\newcommand{\xmath}[1]	{\ensuremath{#1}\xspace}

\newcommand{\reals}	{\xmath{\mathbb{R}}}
\newcommand{\norm}[1]	{\xmath{\left|\left| {#1} \right|\right|}}

\newcommand{\lz}	{\xmath{\ell^0}}
\newcommand{\lo}	{\xmath{\ell^1}}

\newcommand{\fvecsym}	{{\phi}}
\newcommand{\fvec}[1]	{\xmath{\bmath{\fvecsym}_{#1}}}

\newcommand{\fvecmat}	{\xmath{\Phi}}
\newcommand{\fvecmatx}	{\xmath{\widetilde{\Phi}}}
\newcommand{\coefsym}	{{a}}
\newcommand{\coefvec}	{\xmath{\bmath{\coefsym}}}
\newcommand{\coefs}[1]	{\xmath{\coefsym_{#1}}}
\newcommand{\coefvect}[1]	{\xmath{\bmath{\coefsym}(#1)}}
\newcommand{\coefst}[2]	{\xmath{\coefsym_{#1}(#2)}}
\newcommand{\coefvecx}	{\xmath{\bmath{z}}}
\newcommand{\coefsx}[1]	{\xmath{z_{#1}}}
\newcommand{\insig}	{\xmath{\bmath{x}}}

\newcommand{\costf}[1]	{\xmath{C\left(#1\right)}}
\newcommand{\costfj}[1]	{\xmath{\widetilde{C}\left(#1\right)}}
\newcommand{\sigdim}	{\xmath{M}}
\newcommand{\coefdim}	{\xmath{N}}
\newcommand{\tradeoff}	{\xmath{\lambda}}
\newcommand{\sparsity}	{\xmath{S}}
\newcommand{\nodec}	{\xmath{k}}
\newcommand{\nodecb}	{\xmath{j}}
\newcommand{\statesym}	{u}
\newcommand{\statevec}	{\xmath{\bm{\statesym}}}
\newcommand{\states}[1]	{\xmath{\statesym_{#1}}}
\newcommand{\statet}[2]	{\xmath{\statesym_{#1}(#2)}}
\newcommand{\statetder}[2]	{\xmath{\dot{\statesym}_{#1}(#2)}}
\newcommand{\statevect}[1]	{\xmath{\bmath{\statesym}(#1)}}
\newcommand{\statevectder}[1]	{\xmath{\dot{\bmath{\statesym}}(#1)}}
\newcommand{\coefgrp}[1]	{\xmath{\bm{\coefsym}^{#1}}}
\newcommand{\stategrp}[1]	{\xmath{\bm{\statesym}^{#1}}}
\newcommand{\coefsubsym}	{\xmath{\mathcal{A}}}
\newcommand{\coefsub}[1]	{\xmath{\coefsubsym_{#1}}}
\newcommand{\gradcoef}[1]	{\xmath{\nabla_{#1}}}
\newcommand{\subind}	{\xmath{l}}
\newcommand{\tim}	{\xmath{t}}
\newcommand{\timc}	{\xmath{\tau}}
\newcommand{\thresh}	{\xmath{\lambda}}
\newcommand{\threshvec}	{\xmath{\bm{\lambda}}}
\newcommand{\threshs}[1]	{\xmath{\lambda_{#1}}}
\newcommand{\tfunc}[2]	{\xmath{T_{#1}\left(#2\right)}}
\newcommand{\tfuncj}[2]	{\xmath{\widetilde{T}_{#1}\left(#2\right)}}
\newcommand{\tfuncinv}[2]	{\xmath{T^{-1}_{#1}(#2)}}
\newcommand{\tfuncjinv}[2]	{\xmath{\widetilde{T}^{-1}_{#1}\left(#2\right)}}
\newcommand{\energy}	{\xmath{E}}
\newcommand{\energyt}[1]	{\xmath{E(#1)}}
\newcommand{\lbrelax}	{\xmath{\gamma}}

\providecommand{\defstart}{
\begin{tabular}{l|l|l} Macro & Symbol & Meaning \\ \hline }
\providecommand{\defstop}{\end{tabular}}

%% file: abstract.tex

\noindent Recent research has shown that performance in signal processing tasks can often be significantly improved by using signal models based on sparse representations, where a signal is approximated using a small number of elements from a fixed dictionary.   Unfortunately, inference in this model involves solving non-smooth optimization problems that are computationally expensive.  While significant efforts have focused on developing digital algorithms specifically for this problem, these algorithms are inappropriate for many applications because of the time and power requirements necessary to solve large optimization problems.  Based on recent work in computational neuroscience, we explore the potential advantages of continuous time dynamical systems for solving sparse approximation problems if they were implemented in analog VLSI.  Specifically, in the simulated task of recovering synthetic and MRI data acquired via compressive sensing techniques, we show that these systems can potentially perform recovery  at time scales of 10-20$\mu$s, supporting datarates of 50-100 kHz (orders of magnitude faster that digital algorithms).  Furthermore, we show analytically that a wide range of sparse approximation problems can be solved in the same basic architecture, including approximate $\ell^p$ norms, modified \lo norms, re-weighted \lo and $\ell^2$, the block \lo norm and classic Tikhonov regularization.

%% file: intro.tex

\section{Introduction}
\label{sec:intro}

\IEEEPARstart{M}{any} classical approaches to signal and image processing rely on  applying linear filters to incoming data.  This type of processing can be done so efficiently (especially with specialized DSP integrated circuits) that it is  possible to build ``real-time'' systems for many applications.  However, recent research has shown that performance can often be significantly improved by using nonlinear processing
strategies. For example, when presented with imperfect data measurements (e.g., due to noise, blur, missing data, undersampling, etc.), a common approach is to formulate the problem as a regularized inverse problem.  This strategy can be thought of in a Bayesian framework, where the algorithm searches for a signal that was the most likely cause for the measurements, taking into account a prior probability distribution (i.e., a model) on the signal.

While such Bayesian approaches can improve performance in many signal and image processing tasks, these methods rely on solving non-linear optimization problems that are much more computationally expensive than classical linear filtering.  For example, a common family of optimization programs used in this setting minimizes energy functions of the form
\begin{equation}
	\min_{\coefvec} \; \energy = \frac{1}{2}\norm{\insig - \fvecmat \coefvec}^2_2 + \tradeoff\costfj{\coefvec},
\label{eqn:basicopt}	
\end{equation}
where $\insig\in\reals^\sigdim$ is the observed measurement vector, $\coefvec\in\reals^\coefdim$ is a vector representing an estimate of the signal (possibly through coefficients in a transform domain such as Fourier or wavelets), $\fvecmat$ is a $\sigdim\times\coefdim$ matrix representing a linear measurement and corruption process, $\costfj{\cdot}$ is a cost function penalizing $\coefvec$  based on its fit with the signal model, and \tradeoff is a parameter denoting the relative tradeoff between the data fidelity term and the cost function.  Solving this optimization program is equivalent to finding the maximum \emph{a posteriori} (MAP) estimate of the original signal under a Gaussian noise model, with the cost function corresponding to the log prior distribution on the signal.   Basic signal models frequently assume independence among the elements of $\coefvec$, resulting in a cost function that separates into a sum of individual costs $\left(\mbox{i.e., } \costfj{\coefvec} = \sum_\nodec \costf{\coefs{\nodec}}\right)$.  One common example is the $\ell^p$ norm, defined as $\costfj{\coefvec} = \|\coefvec \|_p^p = \left(\sum_i{a_i^p}\right)$.

Significant research activity over the last two decades has focused on signal models based on sparse representations.  In these models, the cost function \costfj{\cdot} is chosen to penalize signals depending on the number of non-zero elements (i.e., the size of the support set of \coefvec).  Sparse representations have drawn significant interest because many natural and man-made signals can be approximated by just a few elements from an appropriately selected basis set~\cite{OLS:1996}.  Because the program in~\eqref{eqn:basicopt} is actually a NP-hard problem when the cost function simply counts the number of non-zero coefficients~\cite{NAT:1995}, much of the recent research has focused either on developing heuristic (often greedy) approximate solutions~\cite{TRO:2004a}, or providing performance guarantees for relaxed versions of the problem~\cite{TRO:2006}.  To date, the strongest theoretical guarantees involve solving the optimization problem in equation~\eqref{eqn:basicopt} when the cost function is the \lo norm
\begin{equation}
	\min_{\coefvec} \;\frac{1}{2}\norm{\insig - \fvecmat \coefvec}^2_2 + \tradeoff\norm{\coefvec}_1,
\label{eqn:bpdn}	
\end{equation}
where $\norm{\coefvec}_1 = \sum_{i=1}^{\coefdim} |\coefs{i}|$.
This optimization program goes by many different names, including Basis Pursuit De-Noising (BPDN) in the signal processing community~\cite{CHE:1999}.  Surprisingly, in many cases of interest it can be shown that solving BPDN recovers the sparsest solution even through~\eqref{eqn:bpdn} is a tractable convex program~\cite{DON:2005}.  



One example of the utility of BPDN is the recent work in compressed/compressive sensing (CS)~\cite{TAO:2006, BAR:2007, CAN:2008b}.  In brief, the CS results give performance guarantees for inverse problems when the signals are highly undersampled $(\sigdim \ll \coefdim)$ and the signal \coefvec is assumed to be sparse (having only $\sparsity<\sigdim$ non-zeros).  The main CS results essentially show that for certain matrices \fvecmat (generally taken to be random), $\sparsity$-sparse signals can be recovered (up to the noise level) by solving BPDN as long as $\sigdim \sim O\left(\sparsity \log ( \coefdim/\sparsity )\right)$.  These results mean that in situations where measurements are costly, a signal can be undersampled during acquisition in exchange for using more computational resources to recover the signal at a later time.  

Despite the long history of optimization in the field of signal processing (see Mattingley \& Boyd~\cite{BOY:2010} for a detailed discussion), the recent advent of applications that utilize optimization directly to perform signal processing tasks (e.g., CS) highlights a specific need for online optimization solvers that can operate in real time or under power constraints.  To mention two example applications that may specifically benefit from real-time or low-power BPDN solutions (respectively), CS techniques have been proposed for both medical imaging~\cite{LUS:2010} and channel estimation for wireless communications~\cite{HAU:2010}. While we will focus on CS as an example application, sparsity-based models (and the corresponding optimization problems) arise in state-of-the-art solutions to problems in a variety of disciplines, including machine learning and computer vision~\cite{WRI:2010a}, as well as signal restoration (e.g., denoising, deblurring, superresolution, inpainting)~\cite{ELA:2010a}.  

Given the importance of solving problems such as BPDN in state-of-the-art algorithms, recent research has focused on dramatically reducing the time it takes to solve this optimization program.  Sparse approximation is particularly challenging because the cost function in~\eqref{eqn:bpdn}, as well as many other cases of interest, is not a smooth function.  Despite much recent progress in developing both fast general purpose convex optimization algorithms~\cite{BOY:2010} and specialized solvers for~\eqref{eqn:bpdn}, these algorithms are unable to solve moderately-sized BPDN problems fast enough to operate in many real-time applications.  In particular, most algorithms for solving BPDN have storage, time and power requirements that  scale unfavorably with the signal size.

Recent work in computational neuroscience has demonstrated a continuous-time dynamical system where the steady-state response is the solution to the program in~\eqref{eqn:basicopt}, and the architecture of the system is designed to efficiently deal with sparsity-inducing cost functions.  Because the dynamics of this system correspond to basic circuit primitives (e.g., leaky integration, simple thresholding, lateral inhibition, etc.), an analog VLSI implementation has the potential to be significantly faster and more power efficient than digital approaches~\cite{MEA1}.
For example, such an implementation could enable applications where CS techniques are used to acquire signals very quickly \emph{and} the signal is recovered virtually instantaneously and with minimal power, thereby eliminating the typical processing bottlenecks of optimization-based signal processing methods (e.g., signal recovery in CS).


The main goal of this paper is to highlight the potential benefits and wide applicability of  analog architectures for efficiently solving sparsity-based optimization programs.  Specifically, this paper makes two main contributions.  First, we provide extensive simulation comparisons of analog systems and digital algorithms for solving BPDN in the context of CS recovery for synthetic and MRI data.  These examples demonstrate that idealized analog architectures could potentially solve individual optimizations at time scales of of 10-20$\mu$s, supporting datarates of 50-100 kHz (orders of magnitude faster that digital algorithms).   Second, we show that a number of other optimization problems arising in the signal processing and statistics communities can be solved using the same basic architecture, including approximate $\ell^p$ norms for $0 \leq  p\leq 1$, modified \lo norms, re-weighted \lo and $\ell^2$, the block \lo norm and classic Tikhonov regularization.  

%% file: background.tex

\section{Background and related work}
\label{sec:back}

\subsection{Dynamical systems for \lo minimization} 
\label{sub:dynsys}

As mentioned above, recent work in computational neuroscience has shown that dynamical systems can be constructed that provably solve the optimization programs in~\eqref{eqn:basicopt} and are efficient for solving the non-smooth problems of interest in sparse approximation.  These systems, known as locally competitive algorithms (LCAs)~\cite{ROZ:2008}, are comprised of a network of analog nodes being driven by the signal to be approximated.  Each node competes with neighboring nodes for a chance to represent the signal, and the steady-state response represents the solution to the optimization problem.  The LCA is a specific type of Hopfield neural network, which have a long history of being used to solve optimization problems~\cite{HOP:1982}.  We note here that other types of network structures have also been proposed recently to approximately solve sparse approximation problems in other ways~\cite{REH:2007,PER:2004}.

Specifically, the $\nodec^{\mathrm{th}}$ node of the LCA is associated with \fvec{\nodec}, the $\nodec^{th}$ column of \fvecmat.  Without loss of generality, we assume each column has unit norm.  This node is described at a given time \tim by an internal state variable \statet{\nodec}{\tim}.  The coefficients \coefvec are related to the internal states \statevec via an activation function \coefvect{\tim} = \tfuncj{\thresh}{\statevect{\tim}} that is parametrized by \thresh.  These activation functions are often taken to be a type of thresholding function.
In the important special case when the cost function is separable, the output of each node \nodec can be calculated independently of all other nodes by a pointwise activation function $\coefst{\nodec}{\tim} = \tfunc{\thresh}{\statet{\nodec}{\tim}}$.   
Individual nodes are leaky integrators driven by an input proportional to $\langle \fvec{\nodec}, \insig \rangle$, and competition between nodes occurs via lateral connections that allow highly active nodes to suppress nodes with less activity.  The dynamics for node \nodec are given by:
\begin{equation}
	\statetder{\nodec}{\tim} = \frac{1}{\timc}\left[ \langle \insig, \fvec{\nodec} \rangle - \statet{\nodec}{\tim}  - \mathop{\sum_{\nodecb=1}^{\coefdim}}_{\nodecb\neq\nodec} \langle \fvec{\nodec}, \fvec{\nodecb}\rangle  \coefst{\nodecb}{\tim}   \right],
\end{equation}
where \timc is the system time constant.  In vector form, the dynamics for the whole network are given by:
\begin{equation}
	\statevectder{\tim} = \frac{1}{\timc}\left[\fvecmat^t \insig - \statevect{\tim} - \left(\fvecmat^t \fvecmat- I\right)\coefvect{\tim}
	\right].
	\label{eqn:LCAvec}
\end{equation}

In~\cite{ROZ:2008} it was shown that for the energy surface \energy given in~\eqref{eqn:basicopt} with a separable cost function, the path induced by the LCA (using the outputs $\coefst{\nodec}{\tim}$ as the optimization variable) ensures $\frac{d\energyt{t}}{dt}\leq 0$ when the cost function satisfies:
\begin{equation}
\lambda \frac{d\costf{\coefs{\nodec}}}{d\coefs{\nodec}} = \states{\nodec} - \coefs{\nodec} = \states{\nodec} -\tfunc{\thresh}{\states{\nodec}} = \tfuncinv{\thresh}{\coefs{\nodec}} - \coefs{\nodec}.
\label{eqn:costthresh}
\end{equation}
The same arguments also extend to the more general case of non-separable cost functions, ensuring $\frac{d\energyt{t}}{dt}\leq 0$ when
\begin{equation}
\thresh \gradcoef{\coefvec}{\costfj{\coefvec}} = \statevec -\coefvec = \statevec - \tfuncj{\thresh}{\statevec} = \tfuncjinv{\thresh}{\coefvec} - \coefvec \label{eq:CoefRelb}.
\end{equation}
Recent followup work~\cite{AUR:2011} establishes stronger guarantees on the LCA, specifically showing that this system is globally convergent to the minimum of $\energy$ (which may be a local minima if \costf{\cdot} is not convex) and proving that the system  converges exponentially fast with an analytically bounded convergence rate.

The relationship in~\eqref{eqn:costthresh} requires cost functions that are differentiable and activation functions that are invertible.  However, the cost function for BPDN (the \lo norm) is non-smooth at the origin and the most effective sparsity-promoting activation functions will likely have non-invertible thresholding properties.  In these cases, one can start with a smooth cost function that is a relaxed version of the desired cost and calculate the corresponding activation function.  Taking the limit of the relaxation parameter in the activation function yields a formula for \tfunc{\thresh}{\cdot} that can be used to solve the desired problem.  Specifically, in the appendix we use the  log-barrier relaxation~\cite{Boy1} to show that the LCA solves BPDN when the activation function is the well-known soft thresholding function:
\[ \costf{\coefs{\nodec}}=|\coefs{\nodec}| \quad\Longleftrightarrow\quad  \coefs{\nodec} = \tfunc{\thresh}{\states{\nodec}} = 
\begin{cases}
0&|\states{\nodec}|\leq\thresh\\
\states{\nodec}-\thresh\mbox{sign}(\states{\nodec})&|\states{\nodec}|>\thresh
\end{cases}.
 \]
Similarly, the LCA can find a local minima to the non-convex optimization program that minimizes the \lz ``norm'' of the coefficients (i.e., number of non-zeros) by using the hard thresholding activation function~\cite{ROZ:2008}:
\[ 
\costf{\coefs{\nodec}}= I\left(\coefs{\nodec}\neq0\right) \quad\Longleftrightarrow\quad
\coefs{\nodec} = \tfunc{\thresh}{\states{\nodec}} = 
\begin{cases}
0&|\states{\nodec}|\leq\thresh\\
\states{\nodec} &|\states{\nodec}|>\thresh
\end{cases},
 \]
where $I(\cdot)$ is the standard indicator function.




\subsection{Digital algorithms for sparse approximation} 
\label{sub:digalgs}

Recent work has focused significant efforts on developing specialized algorithms for solving BPDN on digital platforms.  Several interior point methods have been proposed in this area, including \lo-magic  \cite{CAN:2005} and l1-ls \cite{BOY:2007}.  Alternatively, the GPSR algorithm~\cite{NOW:2007} employs a gradient projection approach to solving the BPDN problem.  Homotopy (or continuation) methods~\cite{MAL:2005,SAL:2010,GAR:2008a} take an entirely different approach, solving a series of optimization problems for a decreasing sequence of tradeoff parameters \tradeoff and utilizing efficient updates to find these sequential solutions.  To speed up the recovery process for very large signals, additional work has sought  to leverage parallel hardware configurations such as multicore~\cite{BRA:2011} and GPU architectures~\cite{LEE:2008}. Multicore processing makes use of the parallelalizable aspects of the algorithm to divide the total computational burden between the available processing units, incurring larger communication overhead for more processors. GPU-based algorithms mainly utilize the ability to perform matrix calculations substantially faster than standard processors. However, while achieving improvements in solution times, neither of these architectures provide favorable scaling properties and it is unclear if they would be able to provide real-time solutions for significantly sized problems.  Also, neither architecture is appropriate for low-power embedded computing applications.



Among digital algorithms, the family of iterative thresholding methods \cite{BLU:2007,DIA:2007,DAU:2004,WOT:2007,FIG:2003a} is most similar to the LCA. These methods iteratively take gradient-type steps to minimize the cost function \eqref{eqn:basicopt} and apply a thresholding function to enforce the sparsity constraints. A first-order discrete Euler approximation to the continuous-time LCA dynamics illustrates that the fundamental update of this analog system is basically the same as these digital algorithms, with the principal difference being that each step of the LCA has an incremental effect on the current solution (rather than taking a large step as in each iteration of the digital algorithm)~\cite{ROZ:2008}.  Recently, approaches based on linearized Bregman iterations have also been shown to have update steps that have a similar form~\cite{OSH:2008}.




%% file: BPDNperformance.tex

\section{Efficient analog BPDN solutions}
\label{sec:bpdnsol}


In this section, we demonstrate the performance of the analog LCA in simulated CS recovery problems to show the potential benefits of analog optimization architectures.  In the first set of simulations (Sections~\ref{sub:phase_plots} and~\ref{sub:convergence_time}), we use synthetic stylized data to thoroughly explore the solution quality and solution times with (simulated) analog and digital approaches.  In the second set of simulations (Section~\ref{sub:MRI}), we use MRI data to show performance on a large scale problem of practical importance.

\subsection{LCA solution quality} 
\label{sub:phase_plots}

To begin, we  investigate the quality of simulated LCA solutions on CS recovery problems with synthetic data to verify that they are comparable  to standard digital algorithms.  While the LCA system is proven to converge asymptotically to the unique BPDN solution, the approximate solution achieved by any algorithm in finite time can have different characteristics depending on the particular solution path.  
In the general problem setup, the unknown signal $a_0 \in \reals^\coefdim$ is $\sparsity$-sparse and is observed through $\sigdim < \coefdim$ Gaussian random projections, $\insig = \fvecmat \coefvec_0 + \nu$, where $\nu$ is additive Gaussian noise.
Following typical approaches in the CS community, we recover an estimate of $a_0$ by solving BPDN.  We compare the simulated performance of the LCA with the interior-point method l1-ls~\cite{BOY:2007} and the gradient projection method GPSR~\cite{NOW:2007}.  This investigation will address two main questions.  First, are the solutions produced by the simulated LCA as accurate as the digital comparison cases?  Second, what solution times are possible in the simulated LCA?  While there may also be significant advantages in power consumption, this issue is implementation specific and beyond the scope of this work.

The test CS problems can be parameterized by the number of observations $\sigdim$, the size of the original sparse signal $\coefdim$ and the sparsity level $\sparsity$. 
We draw the nonzero coefficients of $\coefvec_0$ using a Gaussian distribution with variance $1$ and we draw the locations from a uniform distribution. The choice of regularization parameter $\thresh$ depends on the variance of the additive noise $\nu$ which is not necessarily known a priori. We have empirically observed that $\thresh = .01 \| \fvecmat^T \insig \|_{\infty}$ gives good performance in this task when the noise variance is $10^{-4}$. 
Additionally, we observe that as with many other algorithms, implementing a continuation method by gradually decreasing \thresh (similar to that used in FPC~\cite{WOT:2007}) also improves convergence time in the LCA.  Specifically, we initialize $\thresh =  \| \fvecmat^T \insig \|_{\infty}$ and allow a multiplicative decay of 0.9 at each iteration of the simulation until \thresh reaches the desired value given above.   
To ensure that the comparison among the algorithms is fair, we use the same stopping criterion for convergence based on the duality gap upper bound proposed in~\cite{BOY:2007}.

\begin{figure*}[t]
	\begin{center}
		\includegraphics[width=6in]{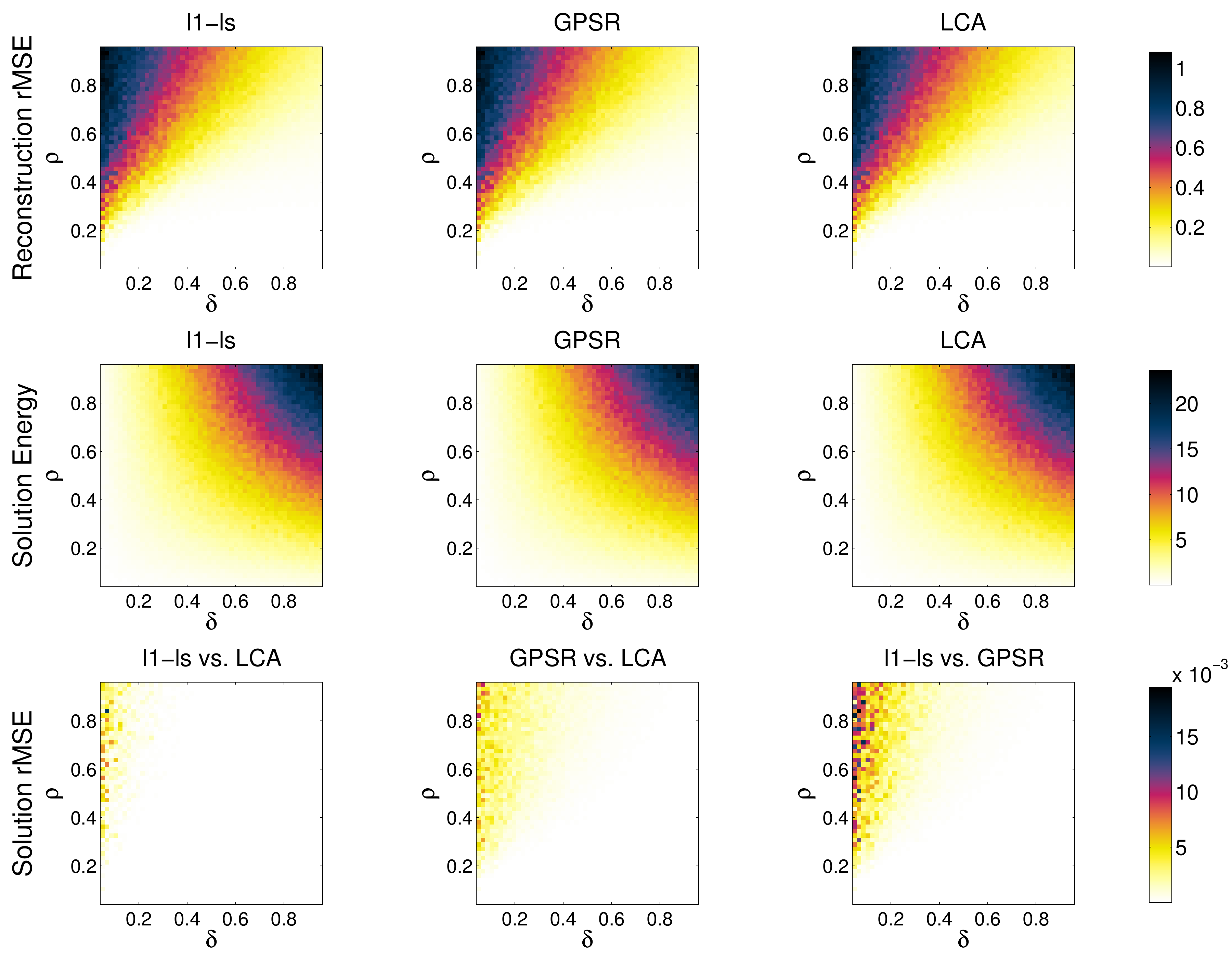} 
	\end{center}
	\caption{The solution quality of the LCA on a compressed sensing recovery task is comparable to the standard digital solvers GPSR and l1-ls. The top row plots the relative MSE of the estimated signal for synthetic data, with indeterminacy of the system indexed by $\delta = \sigdim/\coefdim$, and the sparsity of the system with respect to the number of measurements indexed by $\rho = \sparsity / \sigdim$.  The middle row plots the value of the BPDN objective function at the solutions.  The bottom row plots the relative MSE in the solutions between the solvers, indicating the the differences in the LCA solutions are within the normal range of differences between the digital algorithms themselves.  Note that all solvers demonstrate more variability in regions where the problems are more difficult and signal recovery cannot be performed well. 
}
	\label{fig:EnergyPhase}
\end{figure*}

To explore solution quality we display the results of solving the CS recovery optimizations using plots inspired by the phase plots described by Donoho \& Tanner~\cite{DON:2005}.  We parameterize the plots using the indeterminacy of the system indexed by $\delta = \sigdim/\coefdim$, and the sparsity of the system with respect to the number of measurements indexed by $\rho = \sparsity / \sigdim$. We vary $\delta$ and $\rho$ in the range $[.1, .9]$ using a $50$ by $50$ grid. For a given value $(\delta, \rho)$ on the grid, we sample $10$ different signals using the corresponding $(M, N, S)$ and recover the signal using BPDN.   We compare the results of the simulations by displaying in the top row of Figure~\ref{fig:EnergyPhase} a phase plot for each algorithm, where the color code depicts the average relative MSE of the CS recovery for each algorithm (calculated by $\norm{\hat{\coefvec} - \coefvec_0}_2^2/\norm{\coefvec_0}_2^2$).  In a similar vein, the middle row of Figure~\ref{fig:EnergyPhase} shows the energy function (i.e., the BPDN objective function) evaluated at the solution, $0.5\norm{\insig - \fvecmat\hat{\coefvec}}_2^2 + \thresh\norm{\hat{\coefvec}}_1$.

The near identical plots for the two metrics above demonstrate that the LCA is indeed finding solutions of essentially the same quality as the comparison digital algorithms, both in terms of signal recovery of the compressively sensed signal, and in terms of the optimization objective function.  When the LCA and digital solutions are compared directly, we find that the average difference in the solutions differs only by a relative mean-squared distance (calculated by $\norm{\hat{\coefvec}_{LCA} - \hat{\coefvec}_{DIG}}_2^2/\norm{\hat{\coefvec}_{DIG}}_2^2$) of $1.97\cdot 10^{-4}$ when compared to l1ls and $6.64\cdot 10^{-4}$ when compared to GPSR.  For comparison, the rMSE of the difference between the l1-ls solutions and the GPSR solutions is $9.71\cdot 10^{-4}$, meaning that the LCA solutions have variability comparable to what the pair of comparison digital algorithms has between their solutions.  We note that the solution differences are significantly larger between all of the algorithms in the regimes where CS recovery is difficult and poor solutions are found by all solvers, as demonstrated by the bottom row of plots in Figure~\ref{fig:EnergyPhase}.


%

\subsection{LCA convergence time} 
\label{sub:convergence_time}

\begin{figure*}[t]
	\centering
	\includegraphics[width=5in]{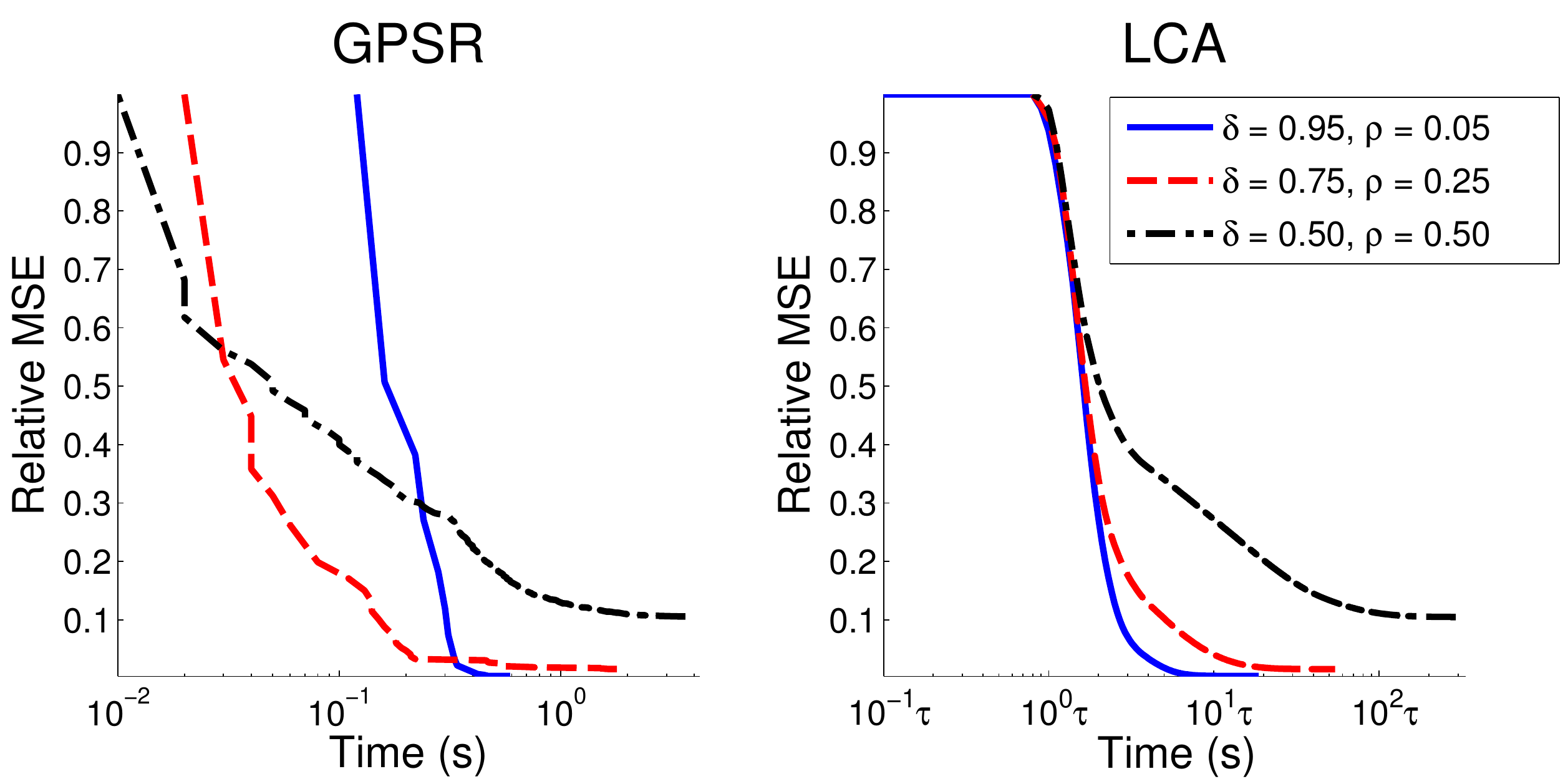}
	\caption{Temporal convergence of the LCA compared to GPSR.  The plot shows the relative MSE of the signal recovery as a function of time for sample trials ($N$=1000) from the results in Figure~\ref{fig:evolution_rmse} using GPSR (left) the simulated LCA (right). The convergence behavior is approximately the same, with harder problems taking both algorithms longer and decreasing the fidelity of the recovery. For the easy and medium difficulty problems where BPDN recovers the signal with good fidelity,  GPSR takes 0.1-1 seconds to converge and the simulated LCA takes $10^1 \tau$-$10^3 \tau$ seconds to converge.  For conservative values of $\tau$, the LCA solution times can still be as low as 10$\mu$s, supporting datarates of up to 100 kHz}
	\label{fig:GPSRandLCAtimes}
\end{figure*}

\begin{figure*}[t]
	\centering
	\includegraphics[width=6in]{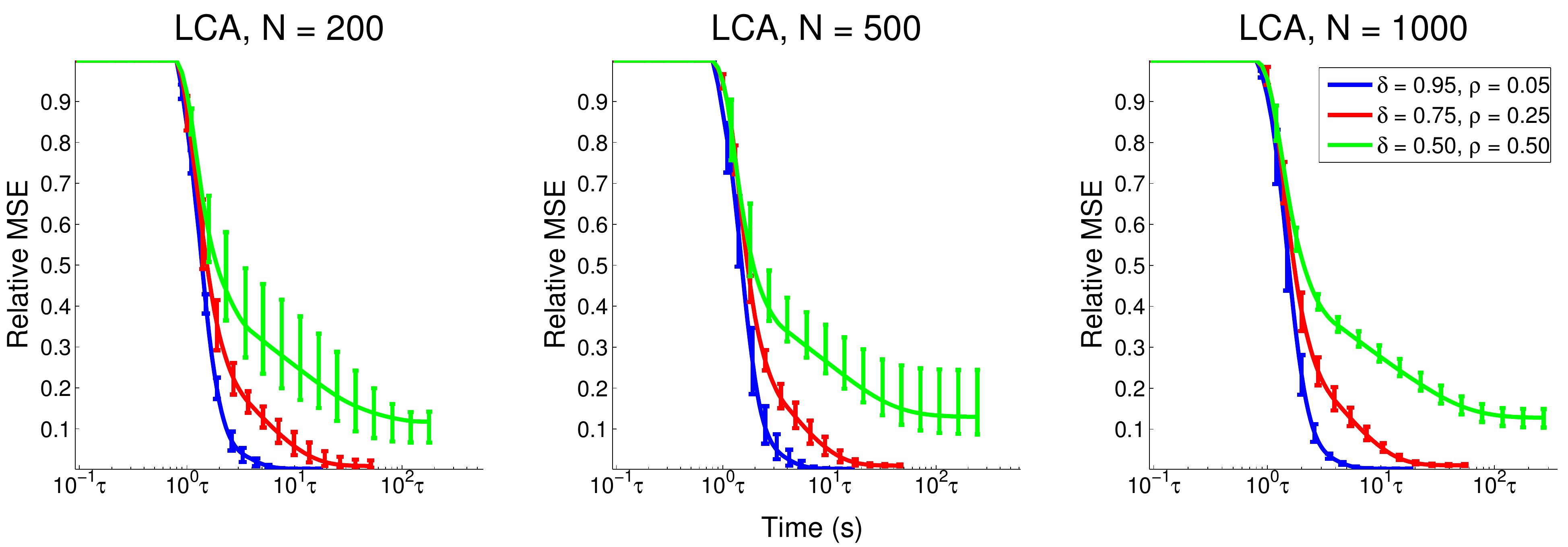}
	\caption{Convergence behavior for the LCA for a number of different problem sizes ($N$,$\delta$,$\rho$). Each plot demonstrates the change in convergence based on easy, medium and hard CS recovery problems (i.e., 3 combinations of ($\delta$, $\rho$)) for $N$ = 200 (left), $N$ = 500 (middle) and $N$ = 1000 (right). While there is no appreciable increase in convergence time with increased problem size (larger $N$), similar to standard behavior with other optimization algorithms the LCA convergence time does increase with problem difficulty (smaller $\delta$ and larger $\rho$).}
	\label{fig:evolution_rmse}
\end{figure*}

To observe the potential solution times for the LCA, we compare the convergence of the LCA and GPSR on three specific signals in easy, medium and hard CS recovery problems with the same synthetic data as above (corresponding to different values of $\delta$, $\rho$).  Figure \ref{fig:GPSRandLCAtimes} shows the convergence of the relative MSE as a function of time for GPSR and the simulated LCA for three example signals.
GPSR times are reported using measured CPU\footnote{Time is measured on a Dell Precision Desktop with dual Intel Xeon E5420 Processors and 14GB of DDR3 RAM.} time, and LCA times are reported using the number of simulated system time constants $\tau$.  The simulation parameters used are identical to the previous simulations.   While the solution paths have generally similar characteristics, the time scales are dramatically different. Focusing on the easy and medium CS problems that produce good recovery using \lo minimization, GPSR is converging in times on the order of 0.3 seconds, whereas the LCA is converging in times on the order of ten time constants ($10\tau$).  We also note that while the results in Figure \ref{fig:GPSRandLCAtimes} are for individual signals for direct comparison with GPSR, the analysis of average case convergence for the LCA shown in Figure~\ref{fig:evolution_rmse} and discussed below also support the same basic conclusions about the LCA convergence time.

Though the time constant of an analog circuit depends on many factors (including the power consumption of the circuit), time constants on the order of 10$^{-6}$ to 10$^{-8}$ are reasonable first-order estimates~\cite{SCH:2011}.  Even with the slowest of these time constants ($\tau \approx 10^{-6}$) the analog solver is converging in approximately 10$\mu$s of simulated time.  This type of solution speed from the LCA 
is several orders of magnitude faster than GPSR and could support solvers running in real time at rates of 100 kHz.  We note that these times are on a similar order as the recent reports of small-scale implementations (especially when accounting for the interface between the analog circuit and the microcontroller hosting the circuit)~\cite{SHA:2011c}.

Finally, we also investigate the effect of problem size $\coefdim$ and problem difficulty ($\delta$, $\rho$) on the convergence speed of the LCA.  For the same parameters corresponding to easy, medium and difficult CS recovery problems as used above, we sample $10$ signals at three different problems sizes ($N$ = 200, $N$ = 500 and $N$ = 1000) to perform CS recovery.  Figure~\ref{fig:evolution_rmse} displays the relative distance of the signal estimate $\coefvec^{(t)}$ from the true solution $\coefvec$ as a function of simulated time, $\norm{\coefvec^{(t)} - \coefvec}_2/\norm{\coefvec}_2$. The plots are again shown as a function of the simulated time in terms of the number of system time constants $\tau$.  As expected, convergence is faster and more reliable (i.e., less variance) for easier recovery problems (i.e., lower sparsity or more measurements).  Interestingly, we note that increasing the signal size \coefdim does \emph{not} appear to increase the solution time for the LCA.  In a digital algorithm such as GPSR, while the number of iterations may not increase substantially, the solution time scales with \coefdim because the cost of each iteration (e.g., a matrix multiplication) increases significantly.  In an analog system like the LCA,  increasing the size of a matrix multiply requires increasing the circuit size and complexity.  While this may increase the system time constant in some implementations~\cite{SHA:2011}, it does not appear to require any more time constants for the system to settle on a solution.\footnote{Note that increasing the problem sizes does increase the time required to simulate the LCA, but not the amount of time being simulated.  Also note that as we will discuss in the conclusions, there may be practical reasons that the system time constant $\tau$ may increase with increasing problem sizes.}

\subsection{MRI Reconstruction}
\label{sub:MRI}

The previous subsection demonstrated that for stylized problems with synthetic data the LCA can achieve BPDN solutions and signal recoveries comparable to standard digital solvers. Furthermore the LCA appears to converge to solutions at speeds that would represent an improvement of several orders of magnitude over digital algorithms if implemented in an analog circuit.  In this section we demonstrate the potential value of this system on a medical imaging application that could be significantly impacted by having real-time CS recovery techniques.  Specifically, in this section we simulate the LCA recovery of undersampled MR images to evaluate the solution quality and speed.  Compressive MRI is of particular interest because it allows shorter scan times, which improves both patient throughput and lowers risk (e.g., shorter scan times mean that pediatric MRIs may be taken more often without general anesthesia~\cite{LUS:2010}).  Furthermore, compressive MR imaging combined with real-time image reconstruction would potentially allow new medical procedures to be performed using real-time, high-resolution 3-D imaging without using ionizing radiation.

We simulate CS data acquisition on 21 frames of a dynamic cardiac MRI sequence\footnote{The MRI data used was acquired using a GE 1.5T TwinSpeed scanner (R12M4) using an 8 element cardiac coil.} by subsampling the Fourier transform of each image (i.e., taking random columns of $k$-space).  Each image is 256x192 pixels, and we recover the images by solving BPDN to find sparse coefficients in a wavelet transform.  Specifically, we solve the BPDN optimization program where the sensing matrix $\bm{\Phi} = \bm{F}\bm{W}^H$ is an inverse wavelet transform followed by a subsampled Fourier matrix, and recover the image by taking the wavelet transform of the solution to the BPDN problem.   The choice of wavelet transforms in this case is very important, as transforms which are coherent with the Fourier subsampling scheme can result in poor results. We follow the work of~\cite{LUS:2010} and use a 4 level 2-dimensional Daubechies wavelet transform as the sparsifying basis.  The resulting optimization is more difficult than the synthetic data in the previous two sections because the signals are larger and the images are sparse in a wavelet basis instead of the canonical basis.

\begin{figure*}[t]
	\centering
	\includegraphics[width=6in]{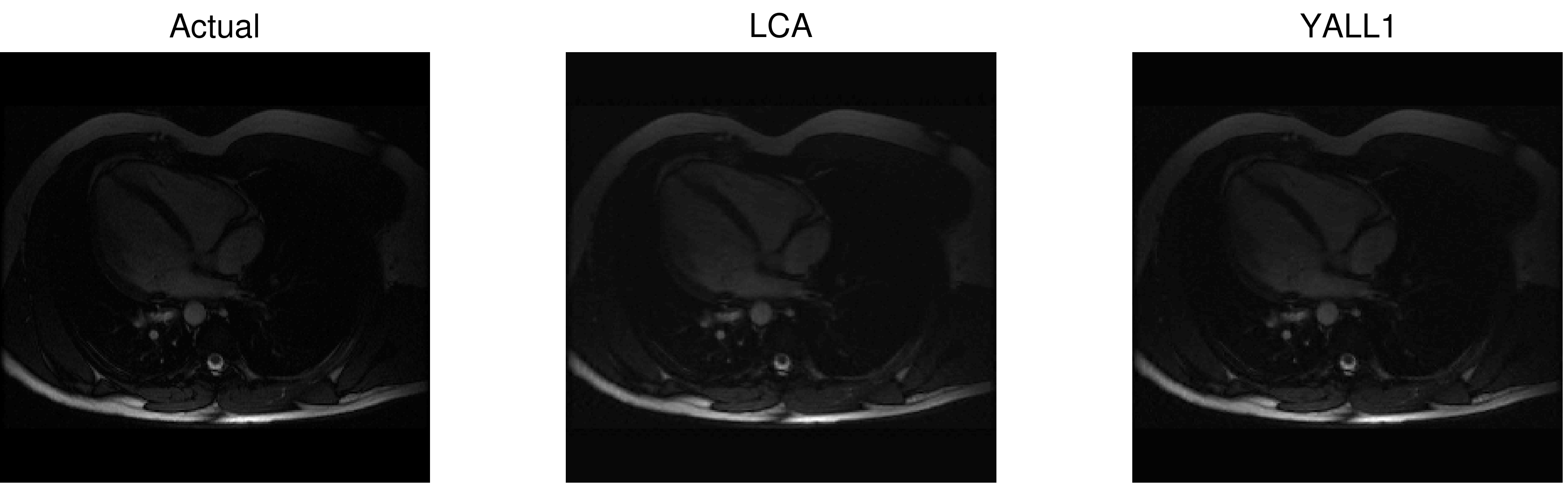}
	\caption{Reconstruction of 256x192 pixel MRI images from simulated CS acquisition.  The simulated LCA and the comparison digital algorithm (YALL1) find solutions of approximately the same quality in terms of relative MSE and image quality. YALL1 finds the solution in approximately 10s, while the LCA finds the solution in approximately 20 time constants ($20 \mu s$ with conservative estimates of the time constant).}
	\label{fig:MRI_recon}
\end{figure*}

We compare results of recovery using the simulated LCA and another standard digital solver YALL1~\cite{WOT:2007}.  Figure~\ref{fig:MRI_recon} shows an example MRI image and its reconstruction using both the LCA and YALL1. The average relative MSE (using $\lambda$ = 0.001) over all 21 recovered images was 0.0109 for YALL1 and 0.0106 for the simulated LCA.  The relative differences between the LCA and YALL1 solutions was 
0.0042, indicating that the solution quality is essentially the same for both approaches. YALL1 took approximately 10 second of computation time to reach this solution (on the same computer platform used in the previous simulations), while the LCA took approximately $20\tau$ simulated seconds.  Again using time constant estimates of $t=10^{-6}$, this translates to solution times of $20 \mu$s and datarates of approximately 50 kHz.

%% file: other_costs.tex

\section{Alternate inference problems in the LCA architecture}
\label{sec:othercosts}

While Section~\ref{sec:bpdnsol} concentrated on exploring the performance of the LCA in solving the commonly used BPDN program, many other cost functions (i.e., signal models) fitting into the general form of~\eqref{eqn:basicopt} have been proposed in the signal processing and statistics literature to exploit sparsity in different ways.  Using the basic relationships described in~\eqref{eqn:costthresh} and~\eqref{eq:CoefRelb}, this section will present a variety of cost functions that can be optimized in the same basic LCA structure by analytically determining the corresponding activation function.\footnote{We also note that a cost function might be easily implementable even in the absence of an analytic formula for the activation function simply by using numerical integration to find a solution and fitting the resulting curve.} These optimization programs include approximate $\ell^p$ norms, modified $\ell^p$ norms that attempt to achieve better statistical properties than BPDN, the group/block \lo norm that induces co-activation structure on the non-zero coefficients, re-weighted \lo and $\ell^2$ algorithms that represent hierarchical statistical models on the coefficients, and classic Tikhonov regularization.

Before exploring specific cost functions, it is worthwhile to make a technical note regarding the optimization programs that are possible to implement in the LCA architecture.  The strong theoretical convergence guarantees established for the LCA~\cite{AUR:2011} apply to a wide variety of possible systems, but do impose some conditions on the permissible activation functions.  We will rely on these same conditions to analytically determine the relationship between the cost and activation functions for the examples we consider in this section.  Translated to conditions on the cost functions, the convergence results for the LCA~\cite{AUR:2011} require that the cost functions be positive $\left(\costfj{\coefvec} \geq 0\right)$, symmetric $\left(\costfj{-\coefvec} = \costfj{\coefvec}\right)$, and satisfy the condition that the matrix $\left(\thresh \gradcoef{\coefvec}^2{\costfj{\coefvec}} + \bm{I}\right)$ is positive definite (i.e., $\thresh \partial^2\costf{\coefs{\nodec}} / \partial\coefs{\nodec}^2 + 1 > 0$ for separable cost functions).   This last condition can intuitively be viewed as requiring that the activation function resulting from~\eqref{eq:CoefRelb} has only a single output for a given input.  In most cases we will only consider the behavior of the activation function for $\states{\nodec} \geq 0$ because the behavior for $\states{\nodec} < 0$ is implied by the symmetry condition.

\begin{figure*}[t]
	\centering
	\includegraphics[width=6in]{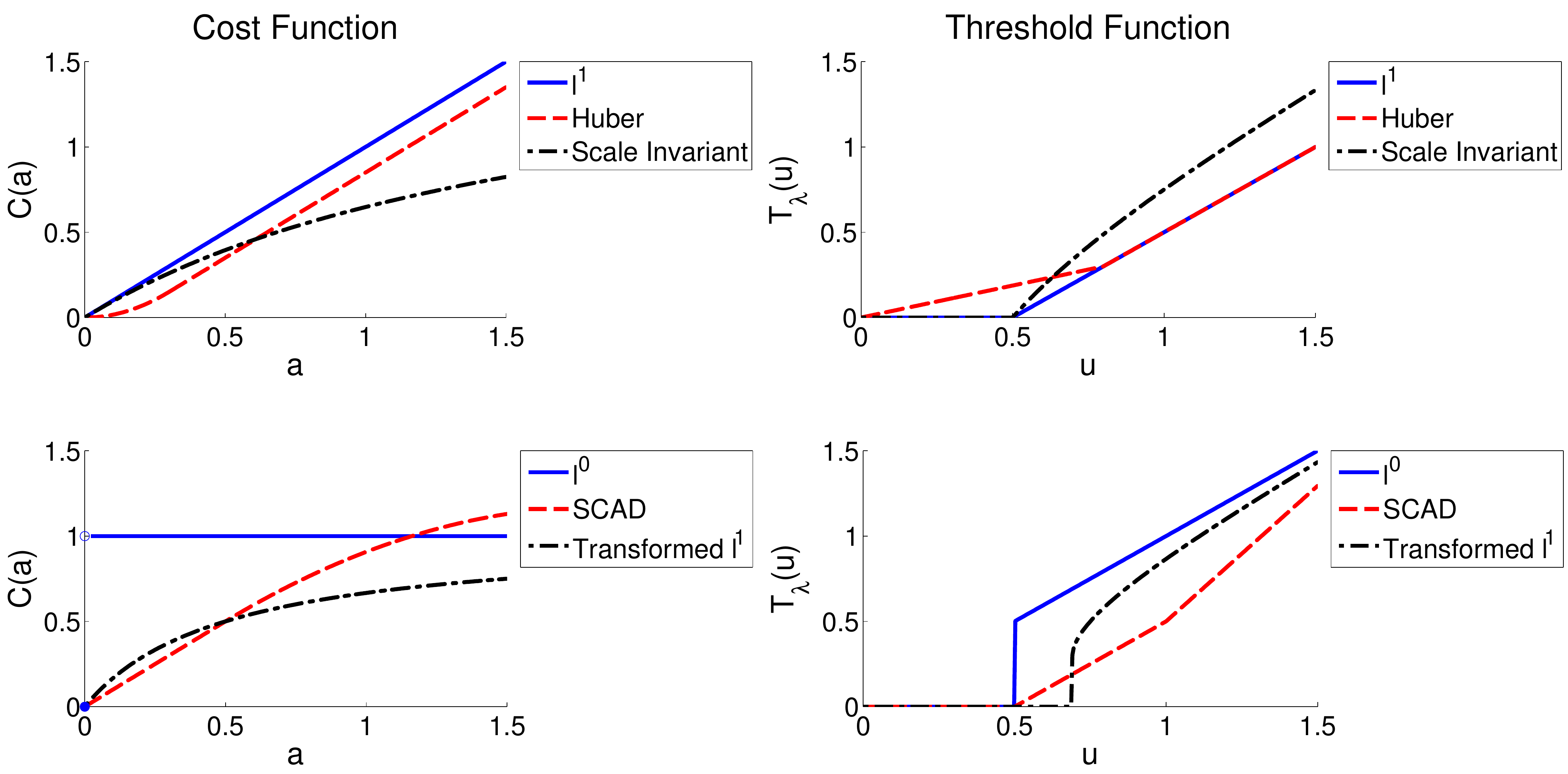}
	\caption{Cost functions and their corresponding thresholding functions. Left: The cost functions are compared for the (top) \lo with \thresh = 0.5, scale invariant Bayes with \thresh = 0.5, the Huber cost with \thresh = 0.5 and $\epsilon$ = 0.3 and (bottom) \lz with \thresh = 0.5, SCAD with \thresh = 0.5 and $\kappa$ = 3.7 and transformed \lo with 
thresh = 0.5 and $\beta$ = 2. Right: The corresponding nonlinear activation function which can be used in the LCA to solve the regularized optimization program for each cost function.}
	\label{fig:costs}
\end{figure*}


\subsection{Approximate $\ell^p$ norms $(0\leq p \leq 2)$} 
\label{sub:approximate_ell_p_norms_0leq_p_leq_2_}

When considering regularized least-squares problems of the form in~\eqref{eqn:basicopt}, perhaps the most widely used family of cost functions are the $\ell^p$ norms $\costfj{\coefvec} = \|\coefvec \|_p^p$.  These separable cost functions include ideal sparse approximation (i.e., counting non-zeros), BPDN, and Tikhonov Regularization~\cite{TIK:1963} as special cases ($p=0,1 \mbox{ and } 2$, respectively), and are convex for $p\geq 1$.  Furthermore, recent research has shown some benefits of using non-convex $\ell^p$ norms ($p<1$) for tasks such as CS recovery~\cite{SAA:2008,ZIB:2007}.  While the ideal activation functions can be determined exactly for the three special cases mentioned above ($p=0,1 \mbox{ and } 2$), it is not possible to analytically determine the activation function for arbitrary values of $0\leq p\leq 2$.  Elad et al.~\cite{ZIB:2007} recently introduced several parameterized approximations to the $\ell^p$ cost functions that are more amenable to analysis.  In this section, we use these same approximations to determine activation functions for minimizing approximate $\ell^p$ norms for $0\leq p\leq 2$.

%
%
%
%
%
%
\subsubsection{Approximate $\ell^p$ for $1\leq p \leq 2$}

\begin{figure*}[t]
	\centering
	\includegraphics[height=2.5in]{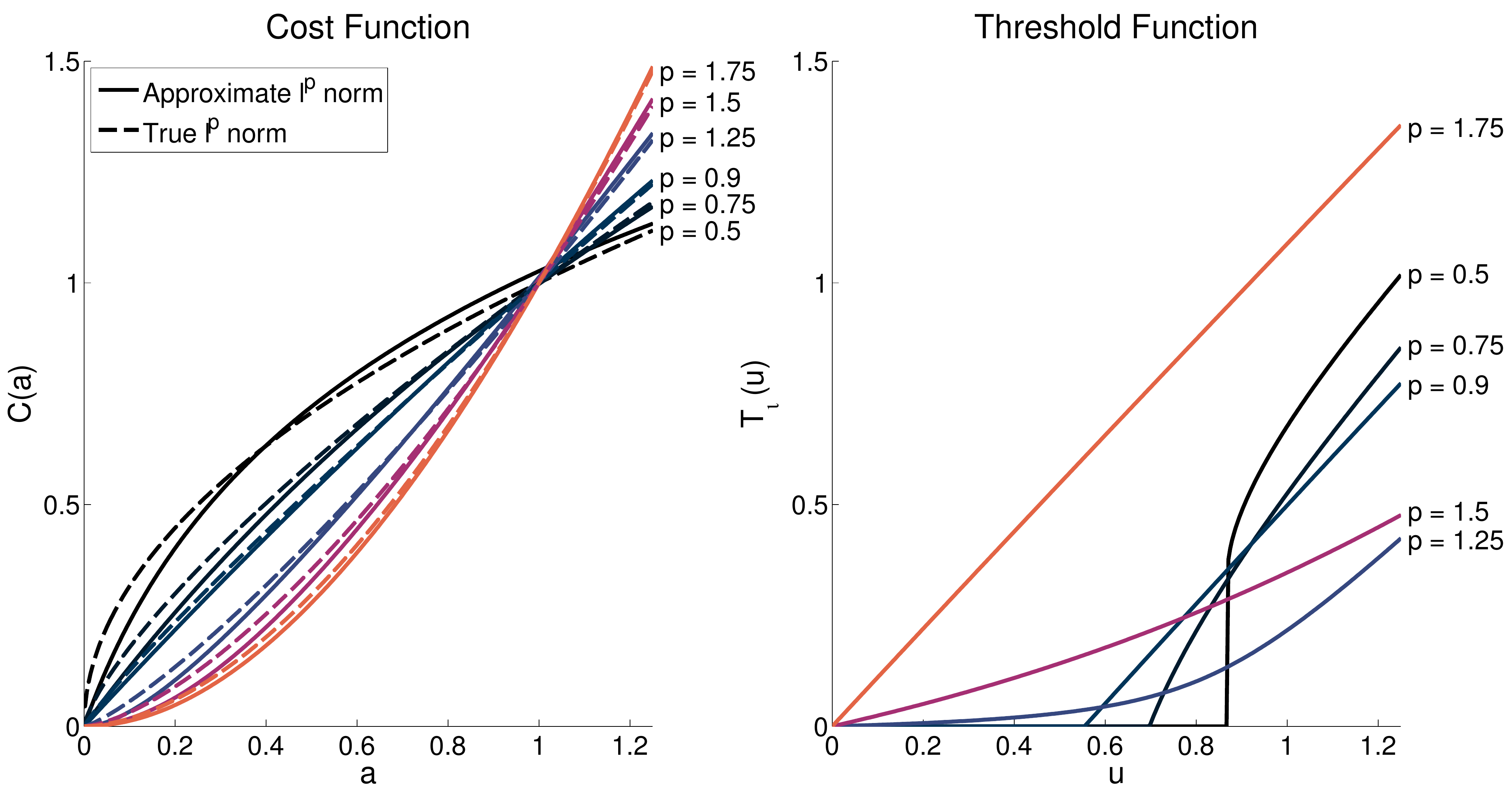}
	\caption{Approximate $\ell^p$ cost functions and their corresponding thresholding functions. Left: The cost functions are approximated over the parameters $c$, $s$ for values of $p$ ranging from 0 to 1 (top) and 1 to 2 (bottom). The true $\ell^p$ costs are shown as dotted lines in the same colors. Using these values of $c$ and $s$, a nonlinear activation function that can be used in the LCA to solve the optimization is plotted (right) using the thresholding equations for $0<p<1$ (top) and $1<p<2$ (bottom). The thresholding functions clearly span the ranges between soft and hard thresholding for the lower range of $p$ and between soft thresholding and linear amplification for the upper range of $p$.}
	\label{fig:costs_lp}
\end{figure*}


For $1 \leq p \leq 2$,  Elad et al.~\cite{ZIB:2007} propose the following approximate cost function as a good match for the true $\ell^p$ norm for some value of parameters $s$ and $c$:
\[ \costf{\coefvec} = \sum_{\nodec}{\left[c|\coefs{\nodec}| - cs\log{\left( 1 + \frac{|\coefs{\nodec}|}{s} \right)}\right]}. \]
In the limiting cases, $c = 1$ with $s \rightarrow 0$ yields the $\lo$ norm and $c = 2s$ with $s \rightarrow \infty$ yields the $\ell^2$ norm. Three intermediate examples for $p$ = 1.25, 1.5 and 1.75 are shown in Figure~\ref{fig:costs_lp}. For any specific value of $p$, we find the best values of $c$ and $s$ by using standard numerical optimization techniques to minimize the squared error to the true cost function over the interval [0,2]. From this cost function, we can differentiate to obtain the relationship between each \states{\nodec} and \coefs{\nodec} as
\[\states{\nodec} = \coefs{\nodec} + \thresh\frac{c\coefs{\nodec}}{s + \coefs{\nodec}}. \]

We see from this relationship that with $c = 1$ and $s \rightarrow 0$, we obtain
$\coefs{\nodec} = \states{\nodec} - \thresh$
for \states{\nodec} $>$ \thresh (i.e., the soft-thresholding function for BPDN), while with $c = 2s$ and $s \rightarrow \infty$ we obtain
$\coefs{\nodec} = \frac{\states{\nodec}}{1 + 2\thresh}$
(i.e., a linear amplifier for Tikhonov Regularization). Solving for $\coefs{\nodec}$ in terms of $\states{\nodec}$ (restricting the solution to be positive and increasing) yields a general relationship for the activation function
\[ \tfunc{\thresh}{\states{\nodec}} = \frac{1}{2}\left[\states{\nodec} -s -c\thresh + \sqrt{\left(\states{\nodec} - s - c\thresh \right) + 4\states{\nodec}s}\right]. \]
This solution is shown in Figure~\ref{fig:costs_lp} for $p$ = 1.25, 1.5 and 1.75 for \thresh = 0.5.


\subsubsection{Approximate $\ell^p$ for $0\leq p \leq 1$}

For $0 \leq p \leq 1$, Elad et al.~\cite{ZIB:2007} also propose the following approximate cost function as a good match for the true $\ell^p$ norm for some value of parameters $s$ and $c$:
\[\costf{\coefs{\nodec}} = cs\log\left(1 + \frac{|\coefs{\nodec}|}{s}\right), \] 
where the parameters $c > 0$ and $s > 0$ can be optimized as above to approximate different values of $p$. Three approximations for $p$ = 0.5, 0.75 and 0.9 are shown in Figure~\ref{fig:costs_lp}. 
To determine the activation function, we again differentiate and find the appropriate relationship to be
%
%
\[\coefs{\nodec} + \frac{\thresh cs}{s+\coefs{\nodec}} = \states{\nodec}. \]
Solving for \coefs{\nodec} reduces to solving a quadratic equation, which leads to two possible solutions. As above, we restrict the activation function to only include the solution that is positive and increasing, resulting in the activation function
\[\tfunc{\thresh}{\states{\nodec}} =  \frac{1}{2}\left(\states{\nodec} - s + \sqrt{\left( \states{\nodec} + s\right)^2 - 4\thresh cs } \right). \]
This activation function is only valid over the range where the output is a positive real number. If $c\thresh\leq s$, this condition reduces to $\states{\nodec}\geq c\thresh$.  More generally, this condition reduces to $\states{\nodec}\geq 2\sqrt{2cs\thresh} - s$.




\subsection{Modified $\ell^p$ norms} 
\label{sub:hybrid_ell_p_norms}

While the general $\ell^p$ norms have historically been very popular cost functions, many people have noted that this approach can have undesirable statistical properties in some instances (e.g., BPDN can result in biased estimates of large coefficients~\cite{ZOU:2006a}).  To  address these issues, many researchers in signal processing and statistics have proposed modified cost functions that attempt to alleviate these statistical concerns.  For example, 
hybrid $\ell^p$ norms smoothly morph between different norms to capture 
the most desirable characteristics over different regions.  In this section we will demonstrate that many of these modified $\ell^p$ norms can also be implemented in the basic LCA architecture.

\subsubsection{Smoothly Clipped Absolute Deviations} 
\label{sub:smoothly_clipped_absolute_deviations}

A common goal for modified $\ell^p$ norms is to retain the continuity of the cost function near the origin demonstrated by the \lo norm, while using a constant cost function for larger coefficients (similar to the \lz norm) to avoid statistical biases.  One approach to achieving these competing goals is the  smoothly clipped absolute deviations (SCAD) penalty~\cite{FAN:1997,FAN:2001}. The SCAD approach directly concatenates the \lo and \lz norms with a quadratic transition region, resulting in the cost function given by
%
%
%
\[\costf{\coefs{\nodec}} = 
\begin{cases}
	\coefs{\nodec} & 0<\coefs{\nodec}\leq\thresh\\
	\frac{1}{(\kappa-1)\thresh}(\coefs{\nodec}\kappa\thresh-\frac{\coefs{\nodec}^2}{2} - \frac{\thresh^2}{2})& \thresh<\coefs{\nodec}\leq\kappa\thresh\\
	\frac{\thresh}{2} (1+\kappa) & \kappa \thresh < \coefs{\nodec}
\end{cases},
\]
for $\kappa \geq 1$ ($\kappa$ defines the width of the transition region).  An example of this cost function with $\thresh = 0.5$ and $\kappa = 3.7$ is shown in Figure~\ref{fig:costs}.

To obtain the activation function we again solve $\thresh \frac{d\costf{\coefs{\nodec}}}{d\coefs{\nodec}} + \coefs{\nodec} = \states{\nodec}$ for \coefs{\nodec} as a function of \states{\nodec}.  For SCAD (and all of the piecewise cost functions we  consider), the activation function can be determined individually for each region, paying careful attention to the ranges of the inputs \states{\nodec} and outputs \coefs{\nodec} to ensure consistency.
%
%
%
For $0<\coefs{\nodec}\leq\thresh$, we have $\thresh + \coefs{\nodec} = \states{\nodec}$, implying  that $\coefs{\nodec}=0$ for $\states{\nodec}<\thresh$ and $\coefs{\nodec}= \states{\nodec}-\thresh$ over the interval $ \thresh<\states{\nodec}<2\thresh$. For $\thresh<\coefs{\nodec}\leq\kappa\thresh$, we have 
\[\thresh\frac{(\kappa\thresh-\coefs{\nodec})}{(\kappa-1)\thresh} + \coefs{\nodec} = \states{\nodec} \implies \coefs{\nodec} = \frac{(\kappa-1) \states{\nodec} - \kappa\thresh }{\kappa-2} \]
over the interval $2\thresh<\states{\nodec}<\kappa\thresh$. Finally, for $\kappa \thresh < \coefs{\nodec}$ we have $ \coefs{\nodec} = \states{\nodec}$, giving the full activation function
\[ \coefs{\nodec} = \tfunc{\thresh}{\states{\nodec}} = 
\begin{cases}
	0 & \states{\nodec}\leq \thresh\\
	\states{\nodec}-\thresh & \thresh \leq \states{\nodec}\leq 2\thresh\\
	\frac{\kappa - 1}{\kappa - 2}\states{\nodec} - \frac{\kappa\thresh}{\kappa - 2} & 2\thresh \leq \states{\nodec} \leq \kappa\thresh\\
	\states{\nodec} & \kappa\thresh \leq \states{\nodec} 
\end{cases},
\]
which is shown in Figure~\ref{fig:costs} for \thresh = 0.5 and $\kappa = 3.7$. Note that this activation function requires $\kappa \geq 2$ (Antoniadis and Fan recommend a value of $\kappa = 3.7$~\cite{FAN:2001}). While this is apparent from consistency arguments once the thresholding function has been derived, this restriction on $\kappa$ can also be deduced from the condition $\thresh\partial^2\costf{\coefs{\nodec}}/\partial\coefs{\nodec}^2 + 1 > 0$.




\subsubsection{Transformed \lo} 
\label{sub:transformed_lo}


Similar to the SCAD cost function, the transformed \lo cost~\cite{FAN:2001,NIK:2000} attempts to capture something close to the \lo norm for small coefficients while reducing the penalty on larger coefficients.  Specifically, transformed \lo uses the fractional cost function given by
\[\costf{\coefs{\nodec}} = \frac{\beta|\coefs{\nodec}|}{1+\beta|\coefs{\nodec}|}, \]
for some $\beta > 0$. An example of this cost with $\beta = 2$ and \thresh = 0.5 is shown in Figure \ref{fig:costs}.  After calculating the derivative of the cost function, the activation function can be found by solving 
%
%
\[ \frac{\thresh\beta}{(1+\beta\coefs{\nodec})^2} + \coefs{\nodec} = \states{\nodec}\]
for \coefs{\nodec}. Inverting this equation reduces to solving a cubic equation in \coefs{\nodec}. The three roots can be calculated analytically, but only one root generates a viable thresholding function by being both positive and increasing for positive \states{\nodec}. That root is given by
\[ \begin{matrix} \coefs{\nodec} = \frac{\beta\, \states{\nodec} - 2}{3\, \beta} + \frac{2^{\frac{2}{3}}}{6\beta} {\left(6\, \beta\, \states{\nodec} - 27\, \beta^2\, \lambda + 6\, \beta^2\, \states{\nodec}^2 + 2\, \beta^3\, \states{\nodec}^3 + 3\, \sqrt{3}\, \beta^3\, \sqrt{-\frac{\lambda\, \left(4\, \beta^3\, \states{\nodec}^3 + 12\, \beta^2\, \states{\nodec}^2 - 27\, \lambda\, \beta^2 + 12\, \beta\, \states{\nodec} + 4\right)}{\beta^4}} + 2\right)}^{\frac{1}{3}} \\
+ \frac{\beta2^{\frac{1}{3}}\, {\left(\beta\, \states{\nodec} + 1\right)}^2}{3\,  {\left(6\, \beta\, \states{\nodec} - 27\, \beta^2\, \lambda + 6\, \beta^2\, \states{\nodec}^2 + 2\, \beta^3\, \states{\nodec}^3 + 3\, \sqrt{3}\, \beta^3\, \sqrt{-\frac{\lambda\, \left(4\, \beta^3\, \states{\nodec}^3 + 12\, \beta^2\, \states{\nodec}^2 - 27\, \lambda\, \beta^2 + 12\, \beta\, \states{\nodec} + 4\right)}{\beta^4}} + 2\right)}^{\frac{1}{3}}} \end{matrix}.\] 
This solution is viable only when \coefs{\nodec} is real valued, which corresponds to the range 
$\states{\nodec} \geq  \left(3\left(\frac{\thresh}{4\beta}\right)^{1/3} - \frac{1}{\beta}\right).$
Outside of this range, no viable non-zero solution exists and so \coefs{\nodec} = 0. The full thresholding function is shown in Figure~\ref{fig:costs} for \thresh = 0.5 and $\beta$ = 2.  While it is interesting that an analytic form can be determined for this activation function, the expression is obviously very complex and would likely have to be approximated by curve fitting in any circuit implementation.

\subsubsection{Huber Function}
\label{sec:hubercost}

The Huber cost function~\cite{HUB:1973} aims to modify standard $\ell^2$ optimization to improve the robustness to outliers. This cost function consists of a quadratic cost function on smaller values and a smooth transition to an \lo cost on larger values, given by
\[\costf{\coefs{\nodec}} = 
\begin{cases}
	\frac{\coefs{\nodec}^2}{2\epsilon} & 0 \leq |\coefs{\nodec}| \leq \epsilon \\
	|\coefs{\nodec}| - \frac{\epsilon}{2} & \epsilon < |\coefs{\nodec}| \\
\end{cases}.
\]
An example of the Huber cost is shown in Figure \ref{fig:costs} for \thresh = 0.5 and $\epsilon$ = 0.3. 
%
As in the case of other piecewise cost functions, we calculate the activation function separately over each interval of interest by calculating the derivative of the cost function in each region. For the first interval, the relationship is given by
$\frac{\thresh\coefs{\nodec}}{\epsilon} = \states{\nodec} - \coefs{\nodec}$, which obviously gives the activation function 
%
$\tfunc{\thresh}{\states{\nodec}} = \frac{\epsilon\states{\nodec}}{\epsilon + \thresh}$ for  $|\states{\nodec}| \leq \epsilon + \thresh$.
For the second interval, we have
$\thresh\frac{\coefs{\nodec}}{|\coefs{\nodec}|} = \states{\nodec} - \coefs{\nodec}$,
which yields the activation function
$\tfunc{\thresh}{\states{\nodec}} = \states{\nodec}\left(1 - \frac{\thresh}{|\states{\nodec}|}\right)$ for $|\states{\nodec}| > \epsilon + \thresh$.  Putting the pieces together, the full activation function (as expected) is a mixture of the Tikhonov regularization and the soft thresholding used for \lo optimization given by
\[\coefs{\nodec} = \tfunc{\thresh}{\states{\nodec}} = 
\begin{cases}
	\frac{\epsilon\states{\nodec}}{\epsilon + \thresh} & |\states{\nodec}| \leq \epsilon + \thresh\\
	\states{\nodec}\left(1 - \frac{\thresh}{|\states{\nodec}|}\right) & |\states{\nodec}| > \epsilon + \thresh\\
\end{cases},
\]
which is shown in Figure~\ref{fig:costs} for \thresh = 0.5 and $\epsilon$ = 0.3. We can see  that as $\epsilon \rightarrow 0$, the cost function converges to the \lo norm and the thresholding function correctly converges back to the soft-threshold function derived earlier using the log-barrier method.


 

\subsubsection{Amplitude Scale Invariant Bayes Estimation} 
\label{sub:ABE}

A known problem with using the \lo norm as a cost function is that it is not scale invariant, meaning that the results can be poor if the amplitude of the input signals changes significantly (assuming a constant value of \thresh).  Many cost functions (including the ones presented above) are heuristically motivated, drawing on intuition and tradeoffs between the behavior of various $\ell^p$ norms.  In contrast, Figueiredo and Nowak~\cite{NOW:2001} approach the problem from the perspective of Bayesian inference with a Jeffreys' prior to determine a cost function with more invariance to amplitude scaling, similar to the non-negative Garrote~\cite{GAO:1998}.  We consider here the cost function 
\[\costf{\coefvec} = \sum_{\nodec} -\frac{\coefs{\nodec}^2}{4\thresh} + \frac{\coefs{\nodec}\sqrt{\coefs{\nodec}^2 + 4\thresh^2}}{4\thresh} + \thresh\log\left( \coefs{\nodec} + \sqrt{\coefs{\nodec}^2 + 4\thresh^2} \right), \]
which is proportional to the one given by Figueiredo and Nowak~\cite{NOW:2001} and is shown in Figure~\ref{fig:costs} for \thresh = 0.5.

Taking the derivative of this cost function, we end up with the relationship between \states{\nodec} and \coefs{\nodec}
\[\states{\nodec} - \coefs{\nodec} = -2\thresh\frac{\coefs{\nodec}}{4\thresh} + \frac{2\thresh}{4\thresh}\sqrt{\coefs{\nodec}^2 + 4\thresh^2}. \]  Solving for \coefs{\nodec} as a function of \states{\nodec} yields the following activation function, 
\[\coefs{\nodec} = \tfunc{\thresh}{\states{\nodec}} = 
\begin{cases}
	0 & \states{\nodec}\leq \thresh\\
	(\states{\nodec}^2-\thresh^2)/\states{\nodec} & \states{\nodec}> \thresh\\
\end{cases}, \]
matching the results from Figueiredo and Nowak~\cite{NOW:2001}.  This activation function is shown in Figure~\ref{fig:costs} for \thresh = 0.5. 
\subsection{Block \lo} 
\label{sub:blockL1}

While all cost functions discussed earlier in this section have been separable,  there is increasing interest in the signal processing community in non-separable cost functions that capture structure (i.e., statistical dependencies) between the non-zero coefficients.  Perhaps the most widely cited cost function discussed in this regard is the block \lo norm (also called the group \lo norm), which assumes that the coefficients representing \insig are active in known groups. In this framework, the coefficients are divided into blocks, $\coefsub{\subind} \subset \left\{ \coefs{\nodec} \right\}$ and each block of coefficients \coefsub{\subind} is represented as a vector \coefgrp{\subind}.  For our purposes, we assume the blocks are non-overlapping but may have different cardinalities. The block \lo norm~\cite{ELD:2010} is defined as the \lo norm over the $\ell^2$ norms of the groups, 
\[ \costfj{\coefvec} = \sum_{\subind}{\left\|\coefgrp{\subind}\right\|_2},	\label{eqn:Cblock} \]
essentially encouraging sparsity between the blocks (i.e., requiring only a few groups to be active) with no individual penalty on the coefficient values within a block.  Because this cost is not separable, the activation function will no longer be a pointwise nonlinearity and will instead have multiple inputs and multiple outputs.

\begin{figure}[t]
	\centering
	\includegraphics[height=2.5in]{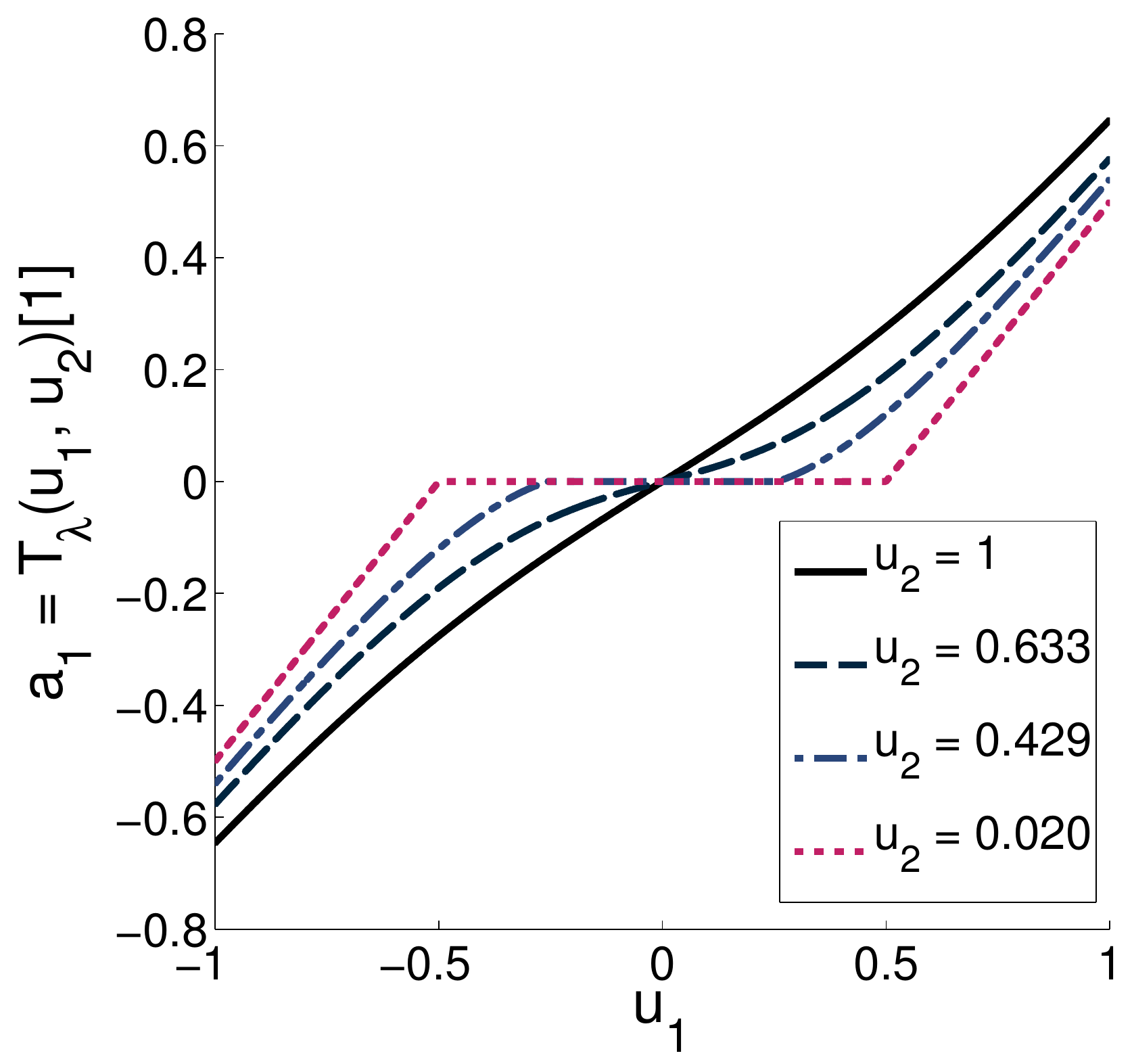} 
	\caption{The nonlinear activation function used in the LCA to optimize the non-overlapping group LASSO cost function has multiple inputs and multiple outputs.  The plot shows an example thresholding function for both elements in a group of size two (\thresh = 0.5), with each line illustrating  the nonlinear effect on \coefs{1} while \states{2} is held constant.}
	\label{fig:GLASSO}
\end{figure}

Following the same general approach as above, we calculate the gradient of the cost function for each block,
\[ \gradcoef{\coefgrp{\subind}}\costfj{\coefvec} = \frac{\coefgrp{\subind}}{\left\|\coefgrp{\subind}\right\|_2}, \]
yielding the following relationship between the activation function inputs and outputs
\begin{equation}
\stategrp{\subind} = \coefgrp{\subind} + \thresh \frac{\coefgrp{\subind}}{\left\|\coefgrp{\subind}\right\|_2}.	
	\label{eqn:group_norms}
\end{equation}
While directly solving this relationship for $\coefgrp{\subind}$ appears difficult, we note that we can simplify the equation by expressing $\left\|\coefgrp{\subind}\right\|_2$ in terms of $\left\|\stategrp{\subind}\right\|_2$.  To see this, take the norm of both sides of~\eqref{eqn:group_norms} to get $\left\|\stategrp{\subind}\right\|_2 = \left\|\coefgrp{\subind}\right\|_2 + \thresh$.  Substituting back into~\eqref{eqn:group_norms}, the relationship simplifies to
%
\[ \tfuncj{\thresh}{\stategrp{\subind}} = \coefgrp{\subind} = \stategrp{\subind}\left(1 - \frac{\thresh}{\left\|\stategrp{\subind}\right\|_2}\right) \]
over the range $0 \leq \left\|\coefgrp{\subind}\right\|_2 = \left\|\stategrp{\subind}\right\|_2 - \thresh$, implying $\thresh \leq \left\|\stategrp{\subind}\right\|_2$.
This relationship yields the block-wise thresholding function
\[\coefgrp{\subind} = \tfuncj{\thresh}{\stategrp{\subind}} = 
\begin{cases}
	0 & \left\|\stategrp{\subind}\right\|_2 \leq \thresh\\
	\stategrp{\subind}\left(1 - \frac{\thresh}{\left\|\stategrp{\subind}\right\|_2}\right) & \left\|\stategrp{\subind}\right\|_2 > \thresh\\
\end{cases}.
\]
This activation function can be thought of as a type of shrinkage operation applied to an entire group of coefficients, with a threshold that depends on the norm of the group inputs.  For the case of groups of two elements (with $\thresh=0.5$),  Figure~\ref{fig:GLASSO} shows the nonlinearities for each of the two states as a function of the value of the other state. 


\subsection{Re-weighted \lo and $\ell^2$} 
\label{sub:reweightedL1}


Recent work has also demonstrated that re-weighted $\ell^p$ norms can achieve better sparsity by iteratively solving a series of tractable convex programs~\cite{WIP:2010,YIN:2008,CAN:2008,GAR:2010}.  For example, re-weighted \lo~\cite{CAN:2008} is an iterative algorithm where a single iteration consists of solving a weighted \lo minimization $\left(\costfj{\coefvec} = \sum_\nodec \thresh_{\nodec} |\coefs{\nodec}|\right)$, followed by a weight update according to the rule
\begin{gather}
	\thresh_\nodec \propto \frac{1}{|\coefs{\nodec}| + \gamma}, \label{eq:rwL1_lambdaupdate}
\end{gather}
where $\gamma$ is a small parameter. By having $\thresh_\nodec$ approximately equal to the inverse of the \lo norm of the coefficient from the previous iteration, this algorithm is more aggressive than BPDN at driving small coefficients to zero and increasing sparsity in the solutions.  
Similarly, re-weighted $\ell^2$ algorithms~\cite{WIP:2010} have also been used to approximate different $p$-norms with weights updated as 
\[ \thresh_\nodec \propto \frac{1}{\left(\coefs{\nodec}^2 + \gamma\right)^{(\frac{p}{2}-1)}}. \] 
Such schemes have shown many empirical benefits over $\ell^p$ norm minimization, and recent work on re-weighted \lo has established theoretical performance guarantees~\cite{KHA:2010} and interpretations as Bayesian inference in a probabilistic model~\cite{GAR:2010}. 

One of the main drawbacks to re-weighted algorithms is the time required for solving the weighted $\ell^p$ program multiple times.  
Because we have established earlier that the LCA architecture can solve the $\ell^p$ norm optimizations (and weighted norms are a straightforward extension to those results), it would immediately follow that a dynamical system could be used to perform the optimization necessary for each iteration of the algorithm.  While this would be a viable strategy (and would save significant time compared to digital solvers, as evidenced by the results in Section~\ref{sub:convergence_time}), we show here that even more advantages can be gained by performing the entire re-weighted \lo algorithm in the context of a dynamical system.  Specifically, we consider here a modified version of the LCA where an additional set of dynamics are placed on \thresh in order to simultaneously optimize the coefficients and coefficient weights in an analog system. While the ideas here are expandable to the general re-weighted case, we focus on results involving the re-weighted \lo as presented in~\cite{GAR:2010}.

\begin{figure}[ht]
	\centering
	\includegraphics[width=6in]{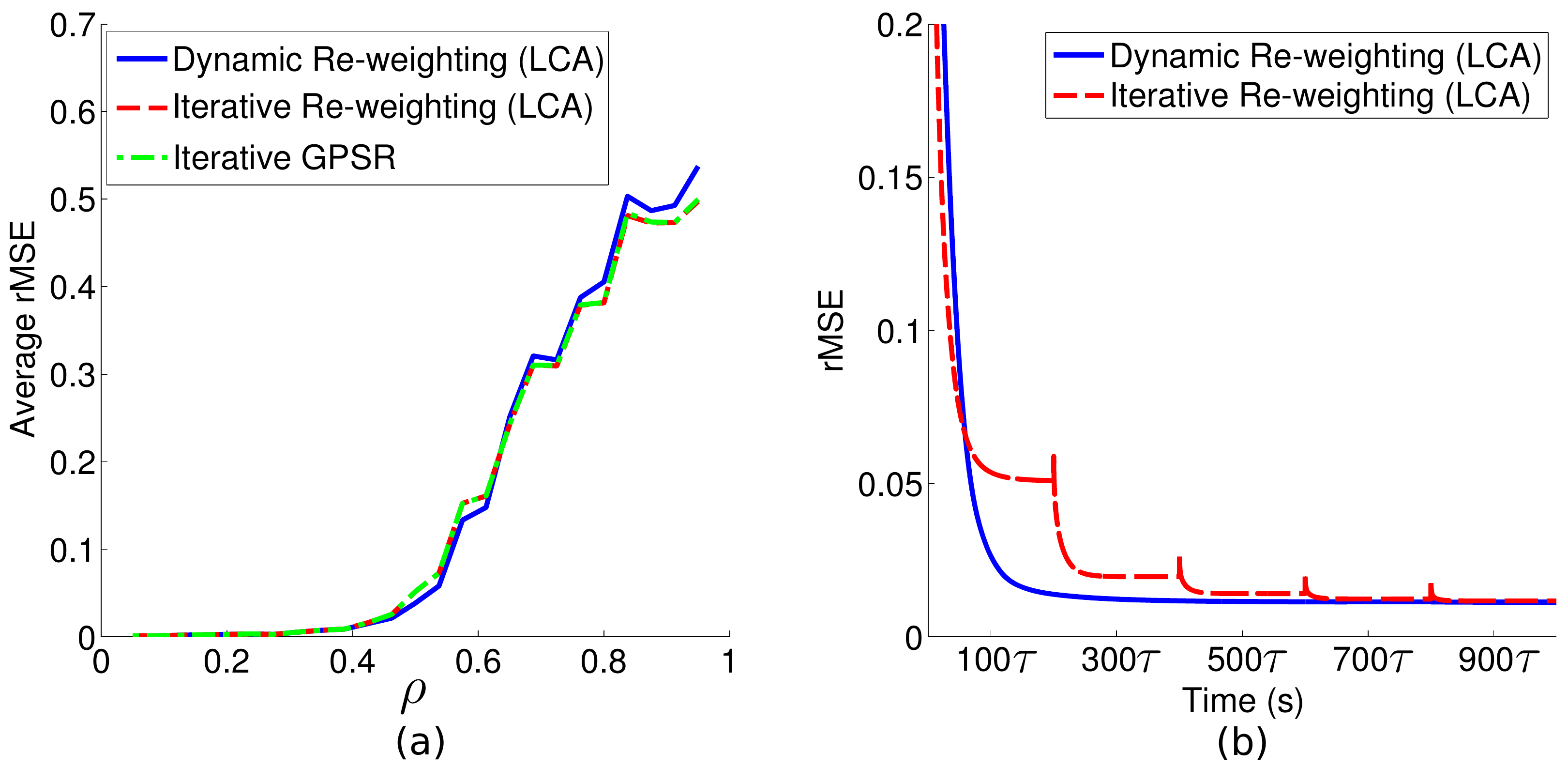}
	\caption{Re-weighted \lo optimization in digital algorithms and in a modified LCA. (a) Re-weighted \lo optimization for a signal with 
	$N = 1000$ and $\delta$ = 0.5, with $\rho$ swept from 0 to 1.  The traditional iterative re-weighting scheme is performed with both a standard digital algorithm (GPSR) and the LCA.  For comparison, a dynamic re-weighting scheme where the LCA is modified to have continuous dynamics on the regularization parameter (rather than discrete iterations) is also shown.  Each method is clearly achieving similar solutions.  (b)  The temporal evolution of the recovery relative MSE for a problem with  $N = 1000$, $\delta$ = 0.6 and $\rho$ = 0.45.  Solutions are shown for the amount of simulated time (in terms of number of time constants).  The dynamically re-weighted system converges in approximately the time it takes to use the LCA to solve two iterations of the traditional  re-weighted \lo algorithm.}
	\label{fig:rwl1mix}
\end{figure}

The modified LCA is given by the system equations:
\[ \begin{matrix} \timc_{\statesym}\dot{\statevec}(t) = \fvecmat^T\insig - \statevec(t) - \left(\fvecmat^T\fvecmat - \bm{I}\right)\coefvec(t) \\ \coefvec(t) = \tfunc{\threshvec}{\statevec(t)} \\ \timc_{\thresh}\dot{\threshs{\nodec}}(t)  = \threshs{\nodec}^{-1}(t) - \nu^{-1}\left(|\coefs{\nodec}(t)| + \gamma\right) \end{matrix}. \]
At steady state, $\dot{\threshvec}$ = 0 which shows that \threshs{\nodec}($\infty$) abides by~\eqref{eq:rwL1_lambdaupdate} with $\nu$ representing the proportionality constant.  While the complete analysis of this expanded analog system is beyond the scope of this paper, we show in Figure~\ref{fig:rwl1mix}a simulations which demonstrate that this system reaches a solution of comparable quality to digital iterative methods. Figure~\ref{fig:rwl1mix}a plots the relative MSE from a CS recovery problem with length-1000 vectors from 500 noisy measurements with varying levels of sparsity.  We sweep the parameter $\rho = S/M$ from zero to one and set the noise variance to $10^{-4}$, with each plot representing the relative MSE averaged over 15 randomly chosen signals.  Figure~\ref{fig:rwl1mix}(a) plots the recovery quality for three systems: iterative re-weighted \lo (using GPSR to solve the \lo iterations), iterative re-weighted \lo (using the LCA to solve the \lo iterations), and dynamic re-weighted \lo which uses the modified LCA described above.  It is clear that the three systems are achieving nearly the same quality in their signal recovery.   Figure~\ref{fig:rwl1mix}b plots the convergence of the recovery as a function of time (in terms of system time constants $\tau$) for the iterative and dynamic re-weighted approaches using the LCA.  The dynamically re-weighted system clearly converges more quickly, achieving its final solution in approximately the time it takes to perform two iterations of the traditional re-weighting scheme using the standard LCA.

%% file: conclusions.tex

\section{Conclusions and future work}
\label{sec:conc}

Sparsity-based signal models have played a central role in many state-of-the-art signal processing algorithms.  The resulting shift toward optimization as a fundamental computational tool in the signal processing toolbox has made it difficult to implement many of these algorithms in applications with significant power constraints or real-time processing requirements.  The main contributions of this paper have been to illustrate the potential advantages of using an analog dynamical system to perform sparse approximation in an analog integrated circuit.  Specifically, our simulations have demonstrated that the idealized LCA could solve problems of significant size on time scales of approximately 10-20$\mu$s, corresponding to real-time solvers at rates approaching 50-100 kHz.  Interestingly, and in stark contrast to using digital algorithms on the same problems, the solution times in the idealized LCA do not appear to scale significantly with the problem size.  Beyond the \lo minimization problem that is most commonly referenced in the literature, we have also demonstrated that the same network structure can implement a wide variety of other cost functions from the signal processing and statistics literature that are related to sparse approximation.

From these results we conclude that solving sparse approximation problems via analog dynamical systems could have a significant impact on a wide range of applications and certainly warrants further investigation.  In the case of CS, the typical mantra has been that CS techniques can help when measurements are expensive and the user is willing to trade reduced measurement burdens for increased computational complexity during signal recovery.  The potential performance of an implementation of the LCA could remove the current bottleneck of CS recovery, making CS techniques applicable in an even wider variety of applications.  With the increased interest in using signal models that incorporate more information than simple signal sparsity (e.g., `structured sparsity' models) for improved CS performance~\cite{BAR:2010}, an interesting avenue for future study would be to develop efficient dynamical systems for performing inference in models with more complex structure than the group \lo norm already established in this paper.

The design and implementation of analog circuits has traditionally been difficult, and it is not immediately clear that the potential benefits of the idealized LCA illustrated in this paper could be achieved in an actual implementation.  As mentioned earlier, the development of reconfigurable analog chips~\cite{TWI:2009} have improved many of the issues related to barriers in the design phase of analog integrated circuits.  In fact, the reconfigurable platform described in~\cite{TWI:2009} has been used to implement a small version of the LCA for solving BPDN~\cite{SHA:2011}.  The preliminary tests of this LCA implementation are on the same order as the simulated solution speeds shown in the present work.

Implementing a system such as the LCA at a scale large enough to be useful in applications will present additional issues that must be addressed in future work.  In particular, the mismatch between elements inherent in the fabrication process and the scaling of the time constant due to factors such as increased load capacitance present challenges that could reduce the effectiveness of the idealized system.  In addition to large scale implementations, interesting future work would include establishing bounds on the solution errors in terms of fabrication mismatch, exploring system designs that exhibit the least potential for time constant increases as the system scales and determining the viability of hybrid analog-digital systems that achieve the benefits of both modalities.  We note here that the initial prototype implementation in~\cite{SHA:2011} reported a system with solutions achieving relative MSE of less than 5\%.

%% file: apx_logbarrier.tex

\appendix[Soft-threshold activation for BPDN using the log-barrier relaxation]
\label{sec:apx_log_barrier}


We will first rewrite the desired BPDN problem in equation~\eqref{eqn:bpdn} in an extended formulation to make the variables non-negative.  Define a new $\sigdim\times2\coefdim$ matrix through the concatenation operation $\fvecmatx=[\fvecmat \;-\fvecmat]$.  Similarly define a vector $\coefvecx=[\coefvecx_+ \; \coefvecx_-]$ of length $2\coefdim$ such that $\coefsx{i}\geq0$ and $\coefvec = \coefvecx_+ - \coefvecx_-$.  Essentially \coefvecx represents the original variables \coefvec by separating them into two subvectors depending on their sign.  We can then write a constrained optimization program that is equivalent to BPDN:
\begin{equation}
\min_{\coefvecx}\frac{1}{2}\norm{\insig-\fvecmatx\coefvecx}_2^2 + \thresh \sum_{\nodec=1}^{2\coefdim} \coefsx{\nodec} \qquad \mbox{s.t.} \quad \coefsx{\nodec}\geq 0.
\label{eqn:exBPDN}
\end{equation}
This reformulation is a standard way to show that \lo cost penalties are equivalent to a linear function in a constrained optimization program.
One can then apply the standard log-barrier relaxation to convert the program in~\eqref{eqn:exBPDN} to an approximately equivalent unconstrained program:
\begin{equation}
\min_{\coefvecx}\frac{1}{2}\norm{\insig-\fvecmatx\coefvecx}_2^2 + \thresh \sum_{\nodec=1}^{2\coefdim} \coefsx{\nodec} + \left(\frac{1}{\lbrelax}\right)\sum_{\nodec=1}^{2\coefdim} \log(\coefsx{\nodec}).
\label{eqn:lbBPDN}
\end{equation}
As $\lbrelax\to\infty$, this program approaches the desired program~\eqref{eqn:exBPDN}.  This relaxation strategy underlies an interior point algorithm (called the barrier method) for solving convex optimization programs, where~\eqref{eqn:lbBPDN} is repeatedly solved with increasing values of \lbrelax~\cite{Boy1}.


\begin{figure*}
	\begin{center}
		\includegraphics[width=5in]{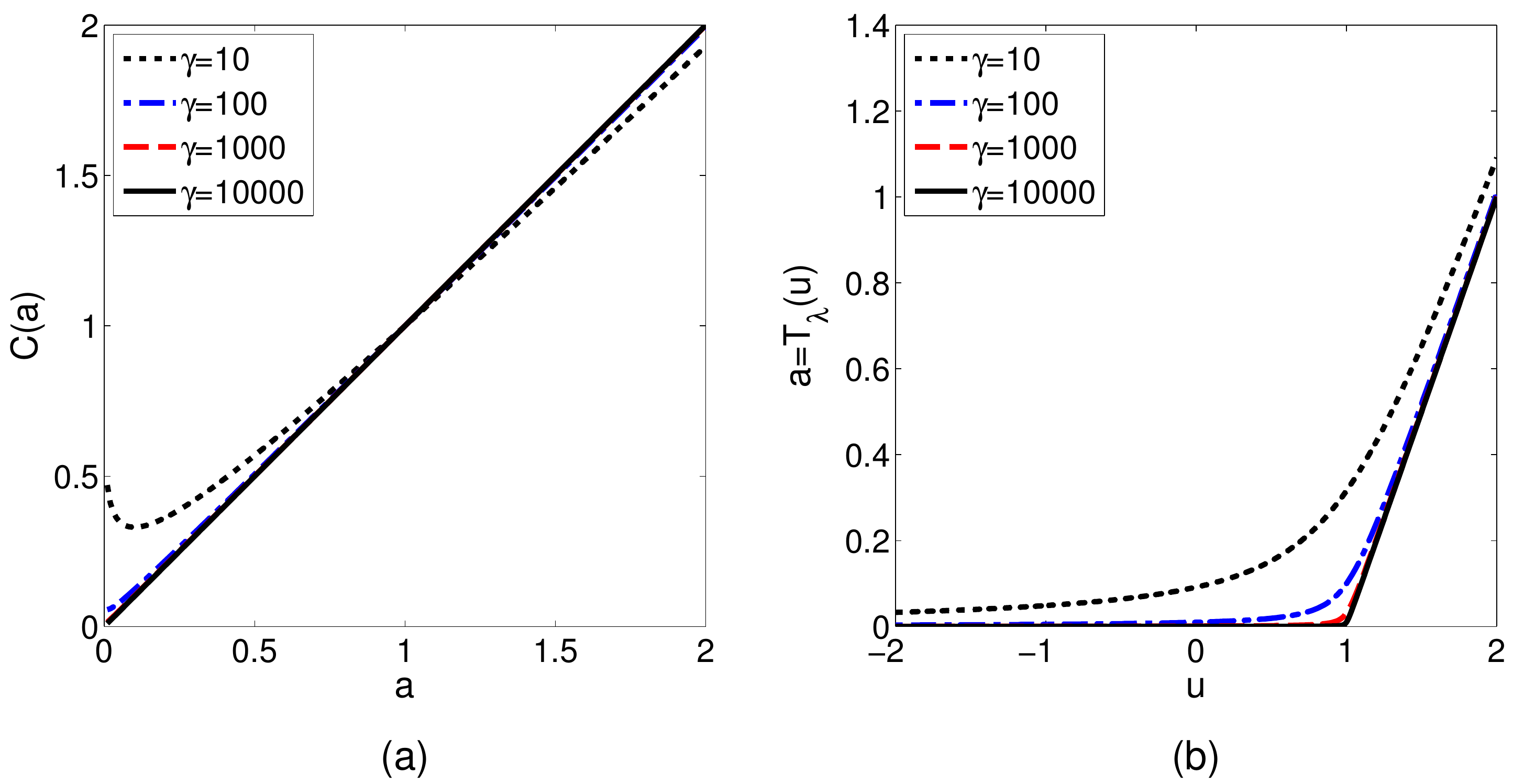}\\
	\end{center}
	\caption{Log barrier relaxations of BPDN.  (a) The cost function approaches the ideal \lo norm as the relaxation parameter is increased.  (b) In a similar way, the nonlinear activation function derived for the LCA approaches the ideal soft-thresholding operator as the relaxation parameter is increased.}
	\label{fig:logbarrier}
\end{figure*}

Note that the relaxed problem in~\eqref{eqn:lbBPDN} fits the form of the general optimization program stated in~\eqref{eqn:basicopt} with the differentiable cost function $\costf{\coefsx{\nodec}} = \coefsx{\nodec}-\frac{\log(\coefsx{\nodec})}{\lbrelax\thresh}.$
For a fixed value of \lbrelax, this cost function can be differentiated and used in the relationship given in~\eqref{eqn:costthresh} to solve for \coefsx{\nodec} in terms of \states{\nodec} to find the corresponding invertible activation function:
\[\coefsx{\nodec} = \tfunc{\thresh}{\states{\nodec}} = \frac{1}{2} \left(\sqrt{\frac{4+\lbrelax (\thresh-\states{\nodec})^2}{\lbrelax}} -(\thresh-\states{\nodec}) \right).\]
Finally it is straightforward to show that in the relaxation limit ($\lbrelax\to\infty$) where the program in~\eqref{eqn:lbBPDN} approaches BPDN, the desired activation function becomes the soft-thresholding function:
\[\lim_{\lbrelax\to\infty} \frac{1}{2} \left(\sqrt{\frac{4+\lbrelax (\thresh-\states{\nodec})^2}{\lbrelax}} -(\thresh-\states{\nodec}) \right) = \frac{1}{2}\left( \sqrt{(\thresh-\states{\nodec})^2} - (\thresh-\states{\nodec})  \right) = 
\begin{cases}
0&\mbox{when }\states{\nodec}\leq\thresh\\
\states{\nodec}-\thresh&\mbox{when }\states{\nodec}>\thresh
\end{cases}.   
\]

To illustrate the convergence of this relaxation to the desired \lo cost function and the corresponding soft-threshold activation function, Figure~\ref{fig:logbarrier} plots \costf{\cdot} and \tfunc{\thresh}{\cdot}
in this relaxed problem for several values of \lbrelax.  Note that in the extended formulation of BPDN given in~\eqref{eqn:exBPDN}, the variables occur in pairs where where only one of them can be nonzero at a time.  Because the activation function is zero for all state values with magnitude less than threshold, it is possible to represent each of these pairs of variables in one LCA node that can take on positive and negative values and where the activation function is a two-sided soft-thresholding function (thereby reducing the number of nodes back down to \coefdim).